\documentclass{amsart}
\usepackage{amssymb,stmaryrd}
\usepackage{amsfonts}
\usepackage{amstext}
\usepackage{algorithmic}
\usepackage{algorithm}
\usepackage{young}
\usepackage{graphicx}
\usepackage{epstopdf}
\usepackage[all]{xy}
\usepackage{enumerate}

\usepackage[colorlinks,linkcolor=blue,citecolor=blue]{hyperref}

\usepackage{color}

\numberwithin{equation}{section}
\newtheorem{theorem}{Theorem}[section]
\newtheorem{lemma}[theorem]{Lemma}
\newtheorem{proposition}[theorem]{Proposition}
\newtheorem{corollary}[theorem]{Corollary}

\theoremstyle{definition}
\newtheorem{definition}[theorem]{Definition}
\newtheorem{example}[theorem]{Example}

\theoremstyle{remark}
\newtheorem{remark}[theorem]{\bf{Remark}}

\newcommand{\R}{{\mathbb{R}}}

\newcommand{\C}{{\mathbb{C}}}

\newcommand{\<}{{\langle}}
\renewcommand{\>}{{\rangle}}

\newcommand{\cg}{{\mathfrak{g}}}
\newcommand{\cm}{{\mathfrak{m}}}

\newcommand{\CYD}{{\mathcal{{\mathcal{M}}}}}

\newcommand{\Lcal}{{\mathcal{L}}}
\newcommand{\Mcal}{{\mathcal{M}}}

\newcommand{\Ad}{{\rm Ad}}
\newcommand{\ad}{{\rm ad}}

\renewcommand{\ker}{{\rm{ker}}}

\newcommand{\tens}{\otimes}
\newcommand{\id}{{\rm id}}
\newcommand{\bo}{{}^{(1)}}
\newcommand{\bt}{{}^{(2)}}

\renewcommand{\o}{{}_{(1)}}
\renewcommand{\t}{{}_{(2)}}
\renewcommand{\th}{{}_{(3)}}

\newcommand{\extd}{{\rm d}}

\newcommand{\und}{\underline}

\newcommand{\la}{{\triangleright}}
\newcommand{\ra}{{\triangleleft}}

\newcommand{\lbiprod}{{>\!\!\!\triangleleft\kern-.33em\cdot}}
\newcommand{\rbiprod}{{\cdot\kern-.33em\triangleright\!\!\!<}}
\newcommand{\lcross}{{>\!\!\!\triangleleft}}
\newcommand{\rcross}{{\triangleright\!\!\!<}}
\newcommand{\rcocross}{{\blacktriangleright\!\!<}}
\newcommand{\lcocross}{{>\!\!\blacktriangleleft}}
\newcommand{\rlbicross}{{\triangleright\!\!\!\blacktriangleleft}}
\newcommand{\lrbicross}{{\blacktriangleright\!\!\!\triangleleft}}
\newcommand{\dcross}{{\bowtie}}

\begin{document}

\title[Differentials on Poisson-Lie groups]{Noncommutative Differentials on Poisson-Lie groups and pre-Lie algebras}


\author{Shahn Majid \& Wen-Qing Tao}
\address{Queen Mary University of London\\
School of Mathematical Sciences, Mile End Rd, London E1 4NS, UK}
\email{s.majid@qmul.ac.uk, w.tao@qmul.ac.uk}
\thanks{The first author was on leave at the Mathematical Institute, Oxford, during 2014 when this work was completed. The second author is supported by the China Scholarship Council}

\begin{abstract}
We show that the quantisation of a connected simply-connected Poisson-Lie group admits a left-covariant noncommutative differential structure at lowest deformation order if and only if the dual of its Lie algebra admits a pre-Lie algebra structure. As an example, we find a pre-Lie algebra structure underlying the standard 3D differential structure on $\C_q[SU_2]$. At the noncommutative geometry level we show that the enveloping algebra $U(\cm)$ of a Lie algebra $\cm$, viewed as quantisation of $\cm^*$, admits a connected differential exterior algebra of classical dimension if and only if $\cm$ admits a pre-Lie algebra. We give an example where $\cm$ is solvable and we extend the construction to the quantisation of tangent and cotangent spaces of Poisson-Lie groups by using bicross-sum and bosonization of Lie bialgebras. As an example, we obtain natural 6D left-covariant differential structures on the bicrossproduct $\C[SU_2]\lrbicross U_\lambda(su_2^*)$. 
\end{abstract}

\maketitle

\section{Introduction}

It is well-known following Drinfeld\cite{Dri} that the semiclassical object underlying quantum groups are Poisson-Lie groups. This means a Lie group together with a Poisson bracket such that the group product is a Poisson map. The infinitesimal notion of a Poisson-Lie group is a Lie bialgebra, meaning a Lie algebra $\cg$ equipped with a `Lie cobracket' $\delta:\cg\to \cg\tens\cg$ forming a Lie 1-cocycle and such that its adjoint is a Lie bracket on $\cg^*$. Of the many ways of thinking about quantum groups, this is a `deformation' point of view in which the coordinate algebra on a group is made noncommutative, with commutator controlled  at lowest order by the Poisson bracket.

In recent years the examples initially provided by quantum groups have led to a significant `quantum groups approach' to noncommutatuve differential geometry in which the next layers of geometry beyond the coordinate algebra are considered, and often classified with the aid of quantum group symmetry. The most important of these is the differential structure, expressed normally as the construction of a bimodule $\Omega^1$ of `1-forms' over the (noncommutative) coordinate algebra and a map $\extd$ for the exterior differential. These are typically extended to a differential graded algebra $(\Omega, \extd)$ of all degrees with $\extd^2=0$. The semiclassical analysis for what this data means at the Poisson level is known to be a {\em Poisson-compatible preconnection} $\nabla$. The systematic analysis in \cite{BegMa:semi} found, in particular, a no-go theorem proving the non-existence of a left and right translation-covariant differential structure of the classical dimension on standard quantum group coordinate algebras $\C_q[G]$ when $G$ is the Lie group of a complex semisimple Lie algebra $\cg$. In \cite{BegMa:twi} was a similar result for the non-existence of ad-covariant differential structures of classical dimension on enveloping algebras of semisimple Lie algebras. Such results tied in with experience at the algebraic level, where one often has to go to higher dimensional $\Omega^1$, and \cite{BegMa:semi,BegMa:twi} also provided an alternative, namely to consider non-associative exterior algebras corresponding to $\nabla$ with curvature. This has been taken up further in \cite{BegMa15}.

The present paper revisits the analysis focusing more clearly on the Lie algebraic structure. For left-covariant differentials we find this time a clean and positive result which both classifies and constructs at semiclassical level differential calculi of the correct classical dimension on quantum groups for which the dual Lie algebra $\cg^*$ has a certain property linked to being solvable. More precisely, our result (Corollary~\ref{preliecorol}) is that the semiclassical data exists if and only if $\cg^*$ admits a pre-Lie structure $\Xi:\cg^*\tens\cg^*\to\cg^*$. Here a pre-Lie structure is a product obeying certain axioms such that the commutator is a Lie algebra, such objects also being called left-symmetric or Vinberg algebras, see \cite{Cartier} and \cite{Bu2} for two reviews. Our result has no contradiction to $\cg$ semisimple and includes quantum groups such as $\C_q[SU_2]$, where we exhibit the pre-Lie structure that corresponds its known but  little understood 3D calculus in \cite{Wor}.  Even better, the duals $\cg^*$ for all quantum groups $\C_q[G]$ are known to be solvable\cite{Ma:mat} and it may be that all solvable Lie algebras admit pre-Lie algebra structures, a question posed by Milnor, see \cite{Bu2}. This suggests for the first time a route to the construction of a left-covariant differential calculus for all $\C_q[G]$, currently an unsolved problem. We build on the initial analysis of this example in \cite{BegMa:semi}. Next, if the calculus is both left and right covariant (bicovariant), we find an additional condition (\ref{Xi-bi}) on $\Xi$ which we relate to infinitesimal or Lie-crossed modules with the coadjoint action, see Theorems~\ref{almostcross} and \ref{zero-curvature}. 

The paper also covers in detail the important case of the enveloping algebra $U(\cm)$ of a Lie algebra $\cm$, viewed as quantisation of $\cm^*$. This is a Hopf algebra so, trivially, a quantum group, and our theory applies with $\cg=\cm^*$ an abelian Poisson-Lie group with its Kirillov-Kostant Poisson bracket. In fact our result in this example turns out to extend canonically to all orders in deformation theory, not just the lowest semiclassical order. We show (Proposition~\ref{envel}) that $U(\cm)$ admits a connected bicovariant differential exterior algebra of classical dimension if and only if $\cm$ admits a pre-Lie algebra. The proof builds on results in \cite{MT}. We do not require ad-invariance but the result again excludes the case that $\cm$ is semisimple since semisimple Lie algebras do not admit pre-Lie structures. The $\cm$ that are allowed do, however, include solvable Lie algebras of the form $[x_i,t]=x_i$ which have been extensively discussed for the structure of `quantum spacetime' (here $x_i, t$ are now viewed as space and time coordinates respectively), most recently in \cite{BegMa14}. In the 2D case we use the known classification of $2$-dimensional pre-Lie structures over $\C$ in \cite{Bu} to classify all possible left-covariant differential structures of classical dimension. This includes the standard calculus previously used in \cite{BegMa14} as well as some other differential calculi in the physics literature\cite{SMK}. We also cover the first steps of the noncommutative Riemannian geometry stemming from the choice of these different differential structures, namely the allowed quantum metrics. The choice of calculus highly constrains the possible quantum metrics, a new phenomenon identified in \cite{BegMa14} and analysed in greater generality in \cite{BegMa15}. The 4D case and its consequences for quantum gravity are explored in our related paper\cite{MaTao:cos}. 

We then apply our theory to the quantisation of the tangent bundle and cotangent bundle of a Poisson-Lie group. In Section 5, we recall the use of the Lie bialgebra $\cg$ of a Poisson-Lie group $G$ to construct the tangent bundle as a bicrossproduct of Poisson-Lie groups and its associated `bicross-sum' of Lie bialgebras\cite{Ma:book}. Our results (see Theorem~\ref{prelie-T}) then suggest a full differential structure,  not only at semiclassical level, on the associated bicrossproduct quantum groups  $\C[G]\lrbicross U_\lambda(\cg^*)$ in \cite{Ma:mat, Ma:bicross, Ma:book}. We prove this in Proposition~\ref{propC(G)U(g*)} and give $\C[SU_2]\lrbicross U_\lambda(su_2^*)$ in detail. Indeed, these bicrossproduct quantum groups were exactly conceived in the 1980s as quantum tangent spaces of Lie groups.  Meanwhile, in Section 6, we use a pre-Lie structure on $\cg^*$ to make $\cg$ into a braided Lie biaglebra\cite{Ma:blie} (Lemma~\ref{blie}). The Lie bialgebra of the cotangent bundle becomes  a `bosonization' in the sense of \cite{Ma:blie} and we construct in some cases a natural preconnection for the semiclassical differential calculus. As before, we cover abelian Lie groups with the Kirillov-Kostant Poisson bracket and quasitriangular Poisson-Lie groups  as examples.

Most of the work in the paper is at the semiclassical level but occasionally we have results about differentials at the Hopf algebra level as in \cite{Wor}, building on our recent work \cite{MT}. In the Appendix, inspired by the construction on cotangent bundle in Section 6, we generalise the construction there to any cocommutative Hopf algebra. This allows us to specialise to the group algebra of a finite group, which is far from the Poisson-Lie setting of the body of the paper. This is a first step in a general braided Hopf algebras construction of differential calculus and exterior algebras to be studied elsewhere.

\section{Preliminaries}
\subsection{Deformation of noncommutative differentials}

We follow the setting in \cite{BegMa:semi}. Let $M$ be a smooth manifold, one can deform commutative multiplication on $C^\infty(M)$ to an associative multiplication $x\bullet y$ where $x\bullet y=xy+O(\lambda).$ Assume that the commutator can be written as $[x,y]_\bullet=x\bullet y-y\bullet x=\lambda\{x,y\}+O(\lambda^2)$ and that we are working in a deformation setting where we can equate order by order in $\lambda$. One can show that $\{\ ,\ \}$ is a Poisson bracket and thus $(M,\{\ ,\ \})$ is a Poisson manifold.

As vector spaces, $n$-forms $\Omega^n(M)$ and exterior algebra $\Omega(M)$ take classical values. But now $\Omega^1(M)$ is assumed to be a bimodule over $(C^\infty (M),\bullet)$ up to order $O(\lambda^2)$ with $x\bullet\tau=x\tau+O(\lambda)$ and $\tau\bullet x=\tau x+O(\lambda)$. 
The deformed $\extd$ operator $\extd _\bullet:\Omega^n(M)\to\Omega^{n+1}(M)$ is given by $\extd_\bullet x=\extd x+O(\lambda)$ and also a graded derivation to order $O(\lambda^2).$
Define $\gamma:C^\infty (M)\tens\Omega^1(M)\to\Omega^1(M)$ by
\begin{equation*}
x\bullet\tau-\tau\bullet x=[x,\tau]_\bullet=\lambda\gamma(x,\tau)+O(\lambda^2).
\end{equation*}

It was shown in~\cite{Haw,BegMa:semi} that being a bimodule up to order $O(\lambda^2)$ requires that map $\gamma$ satisfies
\begin{gather}
\gamma(xy,\tau)=\gamma(x,\tau)y+x\gamma(y,\tau),\label{pre1}\\
\gamma(x,y\tau )=y\gamma(x,\tau)+\{x,y\}\tau.\label{pre2}
\end{gather}
If $\extd_\bullet$ is a derivation up to order $O(\lambda^2),$ then $\gamma$ should also satisfy
\begin{equation}\label{P-C}
	\extd\{x,y\}=\gamma(x,\extd y)-\gamma(y,\extd x).
\end{equation}
\begin{definition}
Any map $\gamma:C^\infty (M)\tens\Omega^1(M)\to\Omega^1(M)$ satisfying (\ref{pre1}) and (\ref{pre2}) is defined to be a \textit{preconnection} on $M$. A preconnection $\gamma$ is said to be \textit{Poisson-compatible} if in addition (\ref{P-C}) is also satisfied, where $\extd: C^\infty(M)\to \Omega^1(M)$ is the usual exterior derivation.
\end{definition}

In view of the properties of a preconnection, one can rewrite ${\nabla}_{\hat{x}}\tau=\gamma(x,\tau)$ for any $\tau\in \Omega^1(M)$, then the map ${\nabla}_{\hat{x}}:\Omega^1(M)\to\Omega^1(M)$ can be thought of as a usual covariant derivative $\nabla_{\hat{x}}$ along Hamiltonian vector fields $\hat{x}=\{x,-\}$ associated to $x\in C^\infty(M)$ and such that 
${\nabla}_{\hat{x}}(y\tau)=y{\nabla}_{\hat{x}}\tau+\hat{x}(y)\tau,$ which are (\ref{pre1}) and (\ref{pre2}). Then a preconnection $\gamma$ can be viewed as a `partially defined connection' on Hamiltonian vector fields.

From the analysis above, we see that a preconnection controls the noncommutativity of functions and $1$-forms, and thus
plays a vital role in deforming a differential graded algebra $\Omega(M)$ at lowest order, parallel to  the Poisson bracket for $C^\infty(M)$ at lowest order.

\subsection{Poisson-Lie groups and Lie bialgebras}

In this paper, we mainly work over a Poisson-Lie group $G$ and its Lie algebra $\cg$. By definition, the Poisson bracket $\{\ ,\ \}:C^\infty(G)\tens C^\infty(G)\to C^\infty(G)$ is determined uniquely by a so-called Poisson bivector $P=\omega\bo\tens\omega\bt$, i.e., $\{f,g\}=P(\extd f,\extd g)=\omega\bo(\extd f)\omega\bt(\extd g).$
Then $\cg$ is a Lie bialgebra with Lie cobracket $\delta:\cg\to\cg\tens\cg$ given by $\delta(x)=\frac{\extd}{\extd t}\omega\bo(g)g^{-1}\tens \omega\bt(g)g^{-1}|_{t=0}$ where $g=\exp tx \in G$ for any $x\in\cg.$ See~\cite{Ma:book} for more details. The map $\delta$ is a Lie 1-cocycle with respect to the adjoint action, and extends to group 1-cocycles $D(g)=(R_{g^{-1}})_*P(g)$ with respect to the left adjoint action and $D^\vee (g)=(L_{g^{-1}})_*P(g)$ with respect to the right adjoint action. Here $D^\vee(g)=\Ad_{g^{-1}}D(g)$ so these are equivalent. We recall for later that the left group cocycle property here is 
\[ D(uv)=D(u)+\Ad_u(D(v)),\quad \forall u,v\in G,\quad D(e)=0.\]
Given $\delta,$ we can recover $D$ as the unique solution to 
\begin{equation*}
\extd D(\tilde{x})(g)=\Ad_g(\delta x),\quad D(e)=0,
\end{equation*}
where $\tilde{x}$ is the left-invariant vector field corresponding to $x\in\cg$.  We then recover the Poisson bracket by  $P(g)=R_{g*}(D(g))$ for all $g\in G$. 

For convenience, we recall the notion of Lie crossed module~\cite{Ma:blie} of a Lie bialgebra.
\begin{definition}
Let $\cg$ be a Lie bialgebra. A \textit{left $\cg$-crossed module} $(V,\la,\alpha)$ is a left $\cg$-module $(V,\la)$ that admits a left $\cg$-coaction $\alpha:V\to\cg\tens V$ such that
\begin{equation*}
\alpha(x\la v)=([x,\ ]\tens\id+\id\tens x\la)\alpha(v)+\delta(x)\la v
\end{equation*}
for any $x\in\cg,v\in V.$

In the case of $\cg$ is finite-dimensional, the notion of a left $\cg$-crossed module is equivalent to a left $\cg$-module $(V,\la)$ that admits a left $\cg^{*op}$-action $\la'$ satisfying
\begin{equation} \label{liecross}
\phi\o\la' v\<\phi\t,x\>+x\o\la v\<\phi,x\t\>=x\la(\phi\la' v)-\phi\la'(x\la v)
\end{equation}
for any $x\in\cg,\,\phi\in\cg^*$ and $v\in V,$ where the left $\cg^{* op}$-action $\la'$ corresponds to the left $\cg$-coaction $\alpha$ above via $\phi\la' v=\<\phi,v\bo\>v\bt$ with $\alpha(v)=v\bo\tens v\bt.$

We call a left $\cg$-module $V$ with linear map $\la':\cg^*\tens V\to V$ (not necessarily an action) such that (\ref{liecross}) a \textit{left almost $\cg$-crossed module}.
\end{definition}

\subsection{Left-covariant preconnections}
It was studied in \cite[Section 3]{BegMa:semi} that left, right covariance and bicovariance of a differential calculus in terms of the preconnection $\gamma.$ Roughly speaking, $\gamma$ is said to be \textit{left-covariant} (\textit{right-covariant}, or \textit{bicovariant}) if the associated differential calculus on $(C^\infty(G),\bullet)$ is left-covariant (right-covariant, or bicovariant) in the usual sense over $(C^\infty (G),\bullet)$ up to $O(\lambda^2).$

To introduce our results, we explain the notations used in \cite{BegMa:semi} and also give a short review of the results on Poisson-Lie groups there. Firstly, there is one-to-one correspondence between 1-forms $\Omega^1(G)$ and $C^\infty(G,\cg^*)$ the set of smooth sections of the trivial $\cg^*$ bundle. For any 1-form $\tau$, define $\tilde{\tau}\in C^\infty (G,\cg^*)$ by letting $\tilde{\tau}_g=L_g^*(\tau_g).$ Conversely,  any $s\in C^{\infty}(G,\cg^*)$ defines an 1-form (denoted by $\hat{s}$) by setting $\hat{s}_g=L_{g^{-1}}^*(s(g)).$ 
The smoothness of vector fields and $1$-forms and this one-to-one relation can also be shown by Maurer-Cartan form. In particular, we know $\extd a\in\Omega^1(G)$, $\widetilde{\extd a}\in C^\infty(G,\cg^*)$ for any $a\in C^{\infty}(G).$ Denote $\widetilde{\extd a}$ by $\hat{L}_a$,   then $\<\hat{L}_a(g),v\>=\<\widetilde{\extd a}(g),v\>=\<L_g^*((\extd a)_g),v\>=\<(\extd a)_g,{(L_g)}_*v\>={(L_g)}_*(v) a,$ which is the directional derivation of $a$ with respect to $v\in\cg$ at $g$. 

Under above notations, a preconnection $\gamma$  can now be rewritten on $\cg^*$-valued functions as $\tilde{\gamma}:C^\infty(G)\times C^\infty(G,\cg^*)\to C^\infty(G,\cg^*)$ by letting $\tilde{\gamma}(a,\tilde{\tau})=\widetilde{\gamma(a,\tau)}.$ 
Note that for any $\phi,\psi\in\cg^*$ and $g\in G,$ there exist $a\in C^\infty(G)$, $s\in C^\infty(G,\cg^*)$ such that $\hat{L}_a(g)=\phi$ and $s(g)=\psi.$ One can define a map $\widetilde{\Xi}:G\times\cg^*\times\cg^*\to\cg^*$ by 
\begin{equation*}
\tilde{\gamma}(a,s)(g)=\{a,s\}(g)+\widetilde{\Xi}(g,\hat{L}_a(g),s(g)).
\end{equation*} 
For brevity,  the notation for the Poisson bracket is extended to include $\cg^*$ valued function by $\{a,s\}=\omega\bo(\extd a)\Lcal_{\omega\bt}(\hat{s}),$ where $\Lcal$ means Lie derivation.  

It was shown in \cite[Proposition 4.5]{BegMa:semi} that a preconnection $\gamma$ or $\tilde{\gamma}$ is left-covariant if and only if $\widetilde{\Xi}(gh,\phi,\psi)=\widetilde{\Xi}(h,\phi,\psi)$ for any $g,h\in G,\,\phi,\psi\in\cg^*.$ 
Hence for left-covariant preconnection $\tilde{\gamma},$ the map $\widetilde{\Xi}$ defines a map $\Xi:\cg^*\tens\cg^*\to\cg^*$ by $\Xi(\phi,\psi)=\widetilde{\Xi}(e,\phi,\psi)$ and  
\begin{equation}\label{preconnection-xi}
\tilde{\gamma}(a,s)(g)=\{a,s\}(g)+\Xi(\hat{L}_a(g),s(g)).
\end{equation}
Conversely, given $\Xi:\cg^*\tens\cg^*\to\cg^*$ and define $\widetilde{\gamma}$ by the formula (\ref{preconnection-xi}), one can show the corresponding $\gamma$ is a left-covariant preconnection.
In addition, a left-covariant preconnection $\tilde{\gamma}$ is Poisson-compatible if and only if 
\begin{equation}\label{compatible}
\Xi(\phi,\psi)-\Xi(\psi,\phi)=[\phi,\psi]_{\cg^*}.
\end{equation} for any $\phi,\psi\in \cg^*$ \cite[Proposition 4.6]{BegMa:semi}. 

Based on these results, we can write down a formula for preconnection $\gamma$ in coordinates. Let $\{e_i\}$ be a basis of $\cg$ and $\{f^i\}$ dual basis of $\cg^*.$ Denote $\{\omega^i\}$ the basis of left--invariant $1$-forms that is dual to $\{\partial_i\}$ the left-invariant vector fields of $G$ generated by $\{e_i\}.$ Then the Maurer-Cartan form is
\[\omega=\sum_i\omega^i e_i\ \in \Omega^1(G,\cg).\] For any $\eta=\sum_i\eta_i\omega^i\in\Omega^1(G)$ with $\eta_i\in C^\infty(G),$ we know $\eta$ corresponds to $\widetilde{\eta}=\sum_i\eta_i f^i\in C^\infty(G,\cg^*).$ On the other hand, any $s=\sum_i s_i f^i\in C^\infty(G,\cg^*)$ with $s_i\in C^\infty(G)$ corresponds to $\widehat{s}=\sum_i s_i\omega^i\in \Omega^1(G).$ In particular, $\widetilde{\extd a}=\sum_i (\partial_i a) f^i.$

For any $a\in C^\infty(G)$ and $\tau=\sum_i\tau_i \omega^i\in \Omega^1(G),$ we know
$\{a,\widetilde{\tau}\}=\sum_i\{a,\tau_i\}f^i$ and 
$
\Xi(\widetilde{\extd a}(g), \widetilde{\tau}(g))=\Xi(\sum_i (\partial_i a)(g) f^i, \sum_j \tau_j(g) f^j)
=\sum_{i,j} (\partial_i a)(g) \tau_j (g) \Xi (f^i,f^j)
=\sum_{i,j,k} (\partial_i a)(g) \tau_j (g) \<\Xi (f^i,f^j), e_k\>f^k,
$
namely
\begin{align*}
\widetilde{\gamma}(a,\widetilde{\tau})&=\sum_k \left(\{a,\tau_k\}+\sum_{i,j}(\partial_i a) \tau_j \<\Xi(f^i,f^j),e_k\>\right) f^k.\nonumber
\end{align*}
If we write $\Xi^{ij}_k=\<\Xi(f^i,f^j),e_k\>$ (or $\Xi(f^i,f^j)=\sum_k \Xi^{ij}_k f^k$) for any $i,j,k,$
then we have
\begin{equation}\label{preconnection-ijk}
\gamma(a,\tau)=\sum_k \left(\{a,\tau_k\}+\sum_{i,j} \Xi^{ij}_k(\partial_i a) \tau_j\right) \omega^k.
\end{equation}
In particular, we have a more handy formula
\begin{equation}\label{preconnection-ijk-1}
\gamma(a,\omega^j)=\sum_{i,k}(\partial_i a)\<\Xi(f^i,f^j),e_k\>\omega^k=\sum_{i,k}\Xi^{ij}_k(\partial_i a)\omega^k,\quad \forall\,j.
\end{equation}

\section{Bicovariant preconnections}

At the Poisson-Lie group level, it was shown in~\cite[Theorem 4.14]{BegMa:semi} that $\tilde{\gamma}$ is bicovariant if and only if 
\begin{equation}\label{G-cross}
\Xi(\phi,\psi)-\Ad^*_{g^{-1}}\Xi(\Ad^*_g\phi,\Ad^*_g\psi)=\phi(g^{-1}\omega\bo(g))\ad^*_{g^{-1}\omega\bt(g)}\psi,
\end{equation}
for  all $g\in G$ and $\phi,\psi\in\cg^*.$ Working at the Lie algebra level of a connected and simply connected Poisson-Lie group we now have a new result:

\begin{theorem}\label{almostcross}
Let $G$ be a connected and simply connected Poisson-Lie group. The left-covariant preconnection $\tilde{\gamma}$ on $G$ determined by $\Xi:\cg^*\tens\cg^*\to\cg^*$ is bicovariant if and only if 
\begin{equation}\label{g-crossv}
\ad^*_x\Xi(\phi,\psi)-\Xi(\ad^*_x\phi,\psi)-\Xi(\phi,\ad^*_x\psi)=\phi(x\o)\ad^*_{x\t}(\psi),
\end{equation}
where $\delta(x)=x\o\tens x\t,$ namely, $(\ad^*,-\Xi)$ makes $\cg^*$ into a left almost $\cg$-crossed module.
Also, the condition~(\ref{g-crossv}) is equivalent to
\begin{equation}\label{bi}
\delta_{\cg^*}\Xi(\phi,\psi)-\Xi(\phi\o,\psi)\tens\phi\t-\Xi(\phi,\psi\o)\tens\psi\t=\psi\o\tens [\phi,\psi\t]_{\cg^*}.
\end{equation}
\end{theorem}
\proof  We first show the `only if' part.
To get corresponding formula on Lie algebra level for (\ref{G-cross}), we substitute $g$ with $\exp tx$ and differentiate at $t=0.$
Notice that $\frac{\extd}{\extd t}\Ad^*(\exp tx)|_{t=0}=\ad^*_x$ and $\Ad^*(\exp tx)|_{t=0}=\id_{\cg^*}.$ We get 
\begin{equation*}
\ad^*_x\Xi(\phi,\psi)-\Xi(\ad^*_x\phi,\psi)-\Xi(\phi,\ad^*_x\psi)=\phi(x\o)\ad^*_{x\t}(\psi),
\end{equation*} as displayed,
where $\delta(x)=x\o\tens x\t=\frac{\extd}{\extd t}g^{-1}P(g)|_{t=0}$ when $g=\exp tx.$
Now denote $\ad^*_x$ by $x\la$ and $-\Xi(\phi,\ )=\phi\la$, the left $\cg^{*op}$-action, then the left hand side of (\ref{g-crossv}) becomes $-x\la(\phi\la\psi)+\phi\o\la\psi\<\phi\t,x\>+\phi\la(x\la\psi),$ while the right hand side is $\phi(x\o)\ad^*_{x\t}(\psi)=-\phi(x\t)\ad^*_{x\o}(\psi)=-x\o\la\psi\<\phi,x\t\>$.
Hence we obtain 
\begin{equation*}
\phi\o\la\psi\<\phi\t,x\>+x\o\la\psi\<\phi,x\t\>=x\la(\phi\la\psi)-\phi\la(x\la\psi)
\end{equation*}
from (\ref{g-crossv}). This means $\cg^*$ is a left almost $\cg$-crossed module under $(\ad^*,-\Xi).$

Conversely,  we can exponentiate $x$ near zero, and solve the ordinary differential Equation (\ref{g-crossv}) near $g=e$. It has a unique solution (\ref{G-cross}) near the identity.  Since the Lie group $G$ is connected and simply connected, one can show (\ref{G-cross}) is valid on the whole group.

Notice that $\ad^*_x\phi=\phi\o \<\phi\t,x\>$ for any $x\in\cg,\,\phi\in\cg^*,$ so the left hand side of (\ref{g-crossv}) becomes
$-\Xi(\phi,\psi)\o\<\Xi(\phi,\psi)\t,x\>-\Xi(\phi\o,\psi)\<\phi\t,x\>-\Xi(\phi,\psi\o)\<\psi\t,x\>,$ while the right hand side of (\ref{g-crossv}) is 
$\phi(x\o)\psi\o \<\psi\t,x\t\>=\psi\o\<[\phi,\psi\t]_{\cg^*},x\>,$ thus (\ref{g-crossv}) is equivalent to (\ref{bi}) by using non-degenerate pairing between $\cg$ and $\cg^*.$
\endproof

\section{Flat preconnections}
As in~\cite{BegMa:semi}, the \textit{curvature} of a preconnection ${\gamma}$ is defined on Hamiltonian vector fields $\hat{x}=\{x,-\}$ by
\begin{equation*}
R(\hat{x},\hat{y})={\nabla}_{\hat{x}}{\nabla}_{\hat{y}}-{\nabla}_{\hat{y}}{\nabla}_{\hat{x}}-{\nabla}_{\widehat{\{x,y\}}},
\end{equation*}
where $\nabla_{\hat{x}}\tau=\gamma(x,\tau)$ for any $\tau\in\Omega^1(G).$ Note here $[\hat{x},\hat{y}]=\widehat{\{x,y\}}.$ 
The curvature of a preconnection reflects the obstruction to the Jacobi identity on any functions $x,y$ and $1$-form $\tau$ up to third order, namely
\[[x,[y,\tau]_\bullet]_\bullet+[y,[\tau,x]_\bullet]_\bullet+[\tau,[x,y]_\bullet]_\bullet=\lambda^2 R(\hat{x},\hat{y})(\tau)+O(\lambda^3).\] 
This is the deformation-theoretic meaning of curvature in this context. 

We say a preconnection is \textit{flat} if its curvature is zero on any Hamiltonian vector fields, namely \begin{equation*}
\gamma(x,\gamma(y,\tau))-\gamma(y,\gamma(x,\tau))-\gamma(\{x,y\},\tau)=0,
\end{equation*} for all $x,y\in C^\infty(G),\,\tau\in\Omega^1(G).$
This is equivalent to
\begin{equation}\label{flat}
\tilde{\gamma}(x,\tilde{\gamma}(y,s))-\tilde{\gamma}(y,\tilde{\gamma}(x,s))-\tilde{\gamma}(\{x,y\},s)=0
\end{equation}
for all $x,y\in C^\infty(G)$ and $s\in C^\infty(G,\cg^*).$ 

\begin{theorem}\label{zero-curvature}
A Poisson-compatible left-covariant preconnection $\gamma$ on a Poisson-Lie group $G$ with Lie algebra $\cg$ is flat if and only if the corresponding map $-\Xi$ is a right $\cg^*$-action (or left $\cg^{*op}$-action) on $\cg^*$, namely
\begin{equation}\label{leftaction}
\Xi([\phi,\psi]_{\cg^*},\zeta)=\Xi(\phi,\Xi(\psi,\zeta))-\Xi(\psi,\Xi(\phi,\zeta)), \quad\forall\ \phi,\psi,\zeta\in\cg^*.
\end{equation}  In this case, $\gamma$ is bicovariant if and only if the left almost $\cg$-crossed module $\cg^*$ given by $(ad^*,-\Xi)$ in Theorem~\ref{almostcross} is actually a left $\cg$-crossed module.
\end{theorem}
\proof 
Let $\gamma$ be a Poisson-compatible left-covariant preconnection on a Poisson-Lie group $G$. Firstly, we can rewrite formula (\ref{flat}) in terms of $\Xi:\cg^*\tens\cg^*\to\cg^*.$
By definition, the three terms in (\ref{flat}) become
\begin{gather*}
\tilde{\gamma}(x,\tilde{\gamma}(y,s))(g)=\{x,\{y,s\}\}(g)+\{x,\Xi(\hat{L}_y(g),s(g))\}+\Xi(\hat{L}_x(g),\{y,s\}(g))\\
+\,\Xi(\hat{L}_x(g),\Xi(\hat{L}_y(g),s(g))),\qquad\qquad\quad\\
\tilde{\gamma}(y,\tilde{\gamma}(x,s))(g)=\{y,\{x,s\}\}(g)+\{y,\Xi(\hat{L}_x(g),s(g))\}+\Xi(\hat{L}_y(g),\{x,s\}(g))\\
+\,\Xi(\hat{L}_y(g),\Xi(\hat{L}_x(g),s(g))),\qquad\qquad\quad\\
\end{gather*} and 
$$\tilde{\gamma}(\{x,y\},s)(g)=\{\{x,y\},s\}(g)+\Xi(\hat{L}_{\{x,y\}}(g),s(g)).$$
Canceling terms involving the Jacobi identity of Poisson bracket, formula (\ref{flat}) becomes
\begin{gather*}
\{x,\Xi(\hat{L}_y(g),s(g))\}+\Xi(\hat{L}_x(g),\{y,s\}(g))+\Xi(\hat{L}_x(g),\Xi(\hat{L}_y(g),s(g)))\\
-\{y,\Xi(\hat{L}_x(g),s(g))\}-\Xi(\hat{L}_y(g),\{x,s\}(g))-\Xi(\hat{L}_y(g),\Xi(\hat{L}_x(g),s(g)))\\
=\Xi(\hat{L}_{\{x,y\}}(g),s(g)).
\end{gather*}

Note that $\gamma$ is Poisson-compatible, this implies
\begin{align*}
\hat{L}_{\{x,y\}}(g)&=\tilde{\gamma}(x,\hat{L}_y)(g)-\tilde{\gamma}(y,\hat{L}_x)(g)\\
&=\{x,\hat{L}_y\}(g)+\Xi(\hat{L}_x(g),\hat{L}_y(g))-\{y,\hat{L}_x\}(g)-\Xi(\hat{L}_y(g),\hat{L}_x(g)).
\end{align*} and
$\{x,\Xi(\hat{L}_y(g),s(g)\}=\Xi(\{x,\hat{L}_y\}(g),s(g))+\Xi(\hat{L}_y(g),\{x,s\}(g))$ by derivation property of $\{x,-\}.$ 

Then (\ref{flat}) is equivalent to
\begin{equation}\label{flat1}
\begin{split}
\Xi(\hat{L}_x(g),\Xi(\hat{L}_y(g),s(g)))-\Xi(\hat{L}_y(g),\Xi(\hat{L}_x(g),s(g)))\\
=\Xi(\Xi(\hat{L}_x(g),\hat{L}_y(g))-\Xi(\hat{L}_y(g),\hat{L}_x(g)),s(g)).
\end{split}
\end{equation}
for all $x,y\in C^\infty(G)$ and $s\in C^\infty(G,\cg^*).$

Now if $\gamma$ is flat, we can evaluate this equation at the identity $e$ of $G,$ and for any $\phi,\psi,\zeta\in\cg^*,$ set $\phi=\hat{L}_x(e),\,\psi=\hat{L}_y(e)$ and $\zeta=s(e)$ for some $x,y\in C^\infty(G),\,s\in C^\infty(G,\cg^*).$ Then (\ref{flat1}) becomes
\begin{equation*}
\Xi(\Xi(\phi,\psi)-\Xi(\psi,\phi),\zeta)=\Xi(\phi,\Xi(\psi,\zeta))-\Xi(\psi,\Xi(\phi,\zeta)).
\end{equation*}
Use compatibility again, we get (\ref{leftaction}) as displayed. This also shows $\Xi$ is a left $\cg^*$-action on itself, or say, $\cg^*$ is a left $\cg^{*op}$-module via $-\Xi.$

Conversely, if $\cg^*$ is a left $\cg^{*op}$-module via $\la:\cg^*\tens\cg^*\to\cg^*$ and such that $-\phi\la\psi+\psi\la\phi=[\phi,\psi]_{\cg^*},$ i.e., (\ref{leftaction}) holds. This implies (\ref{flat1}) for any $x,y\in C^\infty(G),\,s\in C^\infty(G,\cg^*),$ which is equivalent to (\ref{flat}).

The rest of the proof is immediate from Theorem~\ref{almostcross} when the Poisson-Lie group is connected and simply connected.
\endproof

\subsection{Preconnections and pre-Lie algebras}

Now we recall the notion of \textit{left pre-Lie algebra} (also known as \textit{Vinberg algebra, left symmetric algebra}).
An algebra $(A,\circ)$, not necessarily associative, with product $\circ:A\tens A\to A$ is called a \textit{(left) pre-Lie algebra} if the identity
\begin{equation}\label{left-symmetric}
(x\circ y)\circ z-(y\circ x)\circ z=x\circ(y\circ z)-y\circ(x\circ z).
\end{equation}
holds for all $x,y,z\in A.$
From definition, every associative algebra is a pre-Lie algebra and every pre-Lie algebra $(A,\circ)$ admits a Lie algebra structure (denoted by $\cg_A$) with Lie bracket given by
\begin{equation}
[x,y]_{\cg_A}:=x\circ y-y\circ x
\end{equation}
for all $x,y\in A$. The Jacobi identity of $[\ ,\ ]_{\cg_A}$ holds automatically due to (\ref{left-symmetric}). Regarding this, we can rephrase Theorem~\ref{almostcross} and \ref{zero-curvature} as follows.

\begin{corollary}\label{preliecorol}
A connected and simply connected Poisson-Lie group $G$ with Lie algebra $\cg$ admits a Poisson-compatible left-covariant flat preconnection if and only if $(\cg^*,[\ ,\ ]_{\cg^*})$ admits a pre-Lie structure $\Xi$. Moreover, such left-covariant preconnection is bicovariant if and only if $\Xi$ in addition obeys 
\begin{equation}\label{Xi-bi}
\begin{split}
\delta_{\cg^*}\Xi(\phi,\psi)-\Xi(\phi,\psi\o)\tens\psi\t-\psi\o\tens\Xi(\phi,\psi\t)\\
=\Xi(\phi\o,\psi)\tens\phi\t-\psi\o\tens\Xi(\psi\t,\phi)
\end{split}
\end{equation}
for any $\phi,\psi\in\cg^*.$ 
\end{corollary}
\proof  The first part is shown by (\ref{compatible}) and (\ref{leftaction}). For the bicovariant case, the addition condition on $\Xi$ is (\ref{bi}). Using compatibility and rearranging terms, we know that (\ref{bi}) is equivalent to (\ref{Xi-bi}) as displayed.
\endproof

\begin{example}\label{kk} Let $\cm$ be a finite-dimensional Lie algebra and $G=\cm^*$ be an abelian Poisson-Lie group with its Kirillov-Kostant Poisson-Lie group structure $\{x,y\}=[x,y]$ for all $x,y\in\cm\subset C^\infty(\cm^*)$ or $S(\cm)$ in an algebraic context. By Corollary~\ref{preliecorol}, this admits a Poisson-compatible left-covariant flat preconnection if and only if $\cm$ admits a pre-Lie algebra structure $\circ$. This preconnection is always bicovariant as (\ref{Xi-bi}) vanishes when Lie algebra $\cm^*$ is abelian ($\delta_\cm=0$). Then (\ref{preconnection-ijk}) with $\Xi=\circ$ implies
\[ \nabla_{\hat x}\extd y= \extd (x\circ y),\quad \forall\,x,y\in \cm.\] (Note that $\widetilde{\extd y}$ is a constant-valued function in $C^\infty(G,\cm)$, so $\{x,\widetilde{\extd y}\}\equiv0$ and $\tilde{\gamma}(x,\widetilde{\extd y})=\Xi(x,y)$.)
\end{example}

In fact the algebra and its calculus in this example work to all orders. Thus
the quantisation of $C^\infty(\cm^*)$ is $U_\lambda(\cm)$ regarded as a noncommutative coordinate algebra 
with relations $x y - y x= \lambda [x,y]$. If $\cm$ has an underlying pre-Lie structure then the
above results lead to  relations
\[ [x, \extd y]=\lambda \extd (x\circ y) ,\quad \forall x,y\in \cm\]
and one can check that this works exactly and not only to order $\lambda$ precisely as a consequence
of the pre-Lie algebra axiom. The full result here is:

\begin{proposition}\label{envel}
Let $\cm$ be a finite dimensional Lie algebra over a field $k$ of characteristic zero. Then the enveloping algebra $U(\cm)$ admits a connected bicovariant differential graded algebra with left-invariant $1$-forms $\Lambda^1$ of classical dimension (namely, $\dim \Lambda^1=\dim \cm$) if and only if $\cm$ admits a pre-Lie structure.
\end{proposition}
\proof We refer to \cite{MT,Wor} for the formal definition of a first order differential calculus over a Hopf algebra. A calculus is said to be connected if $\ker\extd =k\{1\}$. It is clear from~\cite[Proposition 2.11, Proposition 4.7]{MT} that such differential graded algebra  on $U(\cm)$ with left-invariant $1$-forms $\cm$ 
correspond to $1$-cocycle $Z_\ra^1(\cm,\cm)$  that extends to a surjective right $\cm$-module map $\omega: U(\cm)^+\to\cm$ and $\ker \extd =k\{1\},$ where the derivation $\extd: U(\cm)\to \Omega^1(U(\cm))=U(\cm)\tens\cm$ is given by $\extd a= a\o\tens\omega(\pi(a\t))$ for any $a\in U(\cm).$
Suppose that $\omega$ is such a map, we take $\zeta=\omega|_\cm\in Z^1_\ra(\cm,\cm).$ For any $x\in\cm$ such that $\zeta(x)=0,$ we have $\extd x=1\tens \omega(x)=0,$ then $\ker\extd=k\{1\}$ implies $x=0,$ so $\zeta$ is an injection hence a bijection as $\cm$ is finite dimensional. Now we can define a product $\circ:\cm\tens\cm\to \cm$ by $x\circ y=-\zeta^{-1}(\zeta(y)\ra x),$ this makes $\cm$ into a left pre-Lie algebra as $[x,y]\circ z=-\zeta^{-1}(\zeta(z)\ra[x,y])=\zeta^{-1}((\zeta(z)\ra y)\ra x)-\zeta^{-1}((\zeta(z)\ra x)\ra y)=x\circ (y\circ z)-y\circ(x\circ z).$

Conversely, if $\cm$ admits a left pre-Lie structure $\circ,$ then $y\ra x=-x\circ y$ makes $\cm$ into a right $\cm$-module and  $\zeta=\id_\cm$ the identity map becomes a bijective $1$-cocycle in $Z^1_\ra(\cm,\cm).$ The extended map $\omega:U(\cm)^+\to\cm$ and the derivation $\extd: U(\cm)\to U(\cm)\tens\cm$ are given by 
$\omega(x_1 x_2 \cdots x_n)=((x_1\ra x_2)\cdots \ra x_n)$ for any $x_1x_2\cdots x_n\in U(\cm)^+$ and 
\[\extd (x_1 x_2 \cdots x_n)=\sum_{p=0}^{n-1}\sum_{\sigma\in Sh(p,n-p)}x_{\sigma(1)}\cdots x_{\sigma(p)}\tens\omega (x_{\sigma(p+1)}\cdots x_{\sigma(n)})\]
for any $x_1x_2\cdots x_n\in U(\cm)$ respectively. We need to show that $\ker\extd =k\{1\}.$ On one hand, $k\{1\}\subseteq \ker\extd,$ as $\extd(1)=0.$ On the other hand,
denote by $U_n(\cm)$ the subspace of $U(\cm)$ generated by the products $x_1x_2\cdots x_p,$ where $x_1,\dots,x_p\in \cm$ and $p\le n.$ Clearly, $U_0=k\{1\},\, U_1(\cm)=k\{1\}\oplus\cm, U_p(\cm)U_q(\cm)\subseteq U_{p+q}(\cm)$ and thus $\left(U_n(\cm)\right)_{n\ge 0}$ is a filtration of $U(\cm).$ 
In order to show $\ker\extd\subseteq k\{1\}$, it suffices to show that the intersection \[\ker\extd \cap U_n(\cm)=k\{1\} \text{ for any integer } n\ge 0.\] We prove  by induction on $n\ge 0$. It is obvious for $n=0,$ and true for $n=1,$ as for any $v=\sum_i x_i\in \ker\extd\cap\cm,$ $0=\extd v=\sum_i 1\tens \omega(x_i)=\sum_i 1\tens x_i$ implies $v=\sum_i x_i=0.$ Suppose that $\ker\extd\cap U_{n-1}(\cm)=k\{1\}$ for $n\ge 2.$ For any $v\in \ker\extd\cap U_n(\cm),$ without loss of generality, we can write $v=\sum_i x_{i_{1}}x_{i_{2}}\cdots x_{i_{n}}+v',$ where $x_{i_{j}}\in \cm$ 
and $v'$ is an element in $U_{n-1}(\cm).$ We have 
\begin{align*}
\extd v&=\sum_i 1\tens \omega(x_{i_{1}}\cdots x_{i_{n}})+\sum_i\sum_{j=1}^n x_{i_{1}}\cdots x_{i_{(j-1)}}\widehat{x_{i_{j}}} x_{i_{(j+1)}}\cdots x_{i_{n}}\tens x_{i_{j}}\\
&+ \sum_{i}\sum_{r=1}^{n-2}\sum_{\sigma\in Sh(r,n-r)} x_{i_{\sigma(1)}}\cdots x_{i_{\sigma(r)}}\tens \omega (x_{i_{\sigma(r+1)}}\cdots x_{i_{\sigma(n)}})+\extd v'.
\end{align*}
We denote the elements $u_{i_{j}}:= x_{i_{1}}\cdots x_{i_{(j-1)}}\widehat{x_{i_{j}}} x_{i_{(j+1)}}\cdots x_{i_{n}}\in U_{n-1}(\cm)$ for any $i, 1\le j\le n.$ Except the term $\sum_i\sum_{j=1}^n u_{i_{j}}\tens x_{i_j},$ all the summands in $\extd v$ lie in $U_{n-2}(\cm)\tens\cm,$ thus $\sum_i\sum_{j=1}^n u_{i_{j}}\tens x_{i_j}$ also lies in $ U_{n-2}(\cm)\tens\cm$ as $\extd v=0.$ This implies $\sum_i\sum_{j=1}^n u_{i_{j}}x_{i_j}=v''$ for some element $v''\in U_{n-1}(\cm).$ Rearrange this and add $(n-1)$'s copies of $\sum_i u_{i_n}x_{i_n}=\sum_i x_{i_1}x_{i_2}\cdots x_{i_n}$ on both sides, we get $n\sum_i (x_{i_1}\cdots x_{i_n})=\sum_i\sum_{j=1}^{n-1} x_{i_1}\cdots x_{i_{(j-1)}}[x_{i_j}, x_{i_{j+1}}\cdots x_{i_n}]+v'',$
therefore \[v=\frac{1}{n}\sum_i\sum_{j=1}^{n-1} x_{i_1}\cdots x_{i_{(j-1)}}[x_{i_j}, x_{i_{j+1}}\cdots x_{i_n}]+\frac{1}{n} v''+v'\in U_{n-1}(\cm).\]

So we show that $v$ actually lies in $\ker\extd\cap U_{n-1}(\cm)$ hence $v\in k\{1\}$ by assumption. Therefore we prove that $\ker\extd\cap U_n(\cm)=k\{1\}$ for any $n\ge 0,$ this finishes our proof.
\endproof

To make contact with real classical geometry in the rest of the paper, the standard approach in noncommutative geometry is to work over $\C$ with  complexified differential forms and functions and to remember the `real form' by means of a $*$-involution. We recall that a differential graded algebra over $\C$ is called \textit{$\ast$-DGA} if it is equipped with a conjugate-linear map $\ast:\Omega\to\Omega$ such that $\ast^2=\id,$ $ (\xi\wedge\eta)^*=(-1)^{|\xi||\eta|}\eta^*\wedge \xi^*$ and $\extd(\xi^*)=(\extd\xi)^*$ for any $\xi,\eta\in\Omega.$ Let $\cm$ be a real pre-Lie algebra, i.e., there is a basis $\{e_i\}$ of $\cm$ with real structure coefficients. Then this is also a real form for $\cm$ as a Lie algebra. In this case $e_i^*=e_i$ extends complex-linearly to an involution $*:\cm \to \cm$ which then makes  $\Omega(U_\lambda(\cm))$ as $*$-DGA if  $\lambda^*=-\lambda,$ i.e., if $\lambda$ is imaginary. If we want $\lambda$ real then we should take $e_i^*=-e_i$. 

\begin{example}\label{2D}
Let $\mathfrak{b}$ be the 2-dimensional complex nonabelian Lie algebra defined by $[x,t]=x.$ It admits~\cite{Bu} five families of mutually non-isomorphic pre-Lie algebra structures over $\mathbb{C}$, which are
\begin{align*}
\mathfrak{b}_{1,\alpha}:&\quad t\circ x=-x,\quad t\circ t=\alpha t;\\
\mathfrak{b}_{2,\beta\neq 0}:& \quad x\circ t=\beta x,\quad t\circ x=(\beta-1) x,\quad t\circ t=\beta t;\\
\mathfrak{b}_3:&\quad t\circ x=-x,\quad t\circ t=x-t;\\ 
\mathfrak{b}_4:&\quad x\circ x= t, \quad t\circ x=-x,\quad t\circ t=-2 t;\\
\mathfrak{b}_5:&\quad x\circ t =x,\quad t\circ t=x+t,
\end{align*} where $\alpha,\beta\in \mathbb{C}.$ (Here $\mathfrak{b}_{1,0}\cong\mathfrak{b}_{2,0}$, so let $\beta\neq 0$.)
Thus there are five families of bicovariant differential calculi over $U_\lambda(\mathfrak{b})$:
\begin{align*}
\Omega^1(U_\lambda(\mathfrak{b}_{1,\alpha})):&\quad [t,\extd x]=-\lambda\extd x,\quad [t,\extd t]=\lambda \alpha \extd t;\\
\Omega^1(U_\lambda(\mathfrak{b}_{2,\beta\neq 0})):& \quad [x,\extd t]=\lambda \beta \extd x,\quad [t,\extd x]=\lambda(\beta-1)\extd x,\quad [t,\extd t]=\lambda\beta\extd t;\\
\Omega^1(U_\lambda(\mathfrak{b}_{3})):&\quad [t,\extd x]=-\lambda\extd x,\quad [t,\extd t]=\lambda\extd x-\lambda \extd t;\\ 
\Omega^1(U_\lambda(\mathfrak{b}_{4})):&\quad [x,\extd x]=\lambda\extd t,\quad [t,\extd x]=-\lambda \extd x,\quad [t,\extd t]=-2\lambda \extd t;\\
\Omega^1(U_\lambda(\mathfrak{b}_{5})):&\quad [x,\extd t]=\lambda \extd x,\quad [t,\extd t]=\lambda\extd x+\lambda\extd t.
\end{align*}
All these examples are $*$-DGA's with $x^*=x$ and $t^*=t$ when $\lambda^*=-\lambda$ as $\{x,  t\}$ is a real form of the relevant pre-Lie algebra. We also need for this that $\alpha,\beta$ are real. 
\end{example}

\begin{example}\label{CqSU2} For $q\in\C, q\neq 0,$ we recall  that the Hopf algebra $\C_q[SL_2]$ is, as an algebra, a quotient of free algebra $\C\<a,b,c,d\>$ modulo relations 
\begin{gather*}
ba=qab,\quad ca=qac,\quad db=qbd,\quad dc=qcd,\quad bc=cb,\\
ad-da=(q^{-1}-q) bc,\quad ad-q^{-1}bc=1.
\end{gather*}
Usually, the generators $a,b,c,d$ are written as a single matrix 
$\begin{pmatrix}
a & b\\
c & d
\end{pmatrix}.$ The coproduct, counit and antipode of $\C_q[SL_2]$ are given by
\[\Delta 
\begin{pmatrix}
a & b\\
c & d
\end{pmatrix}
=
\begin{pmatrix}
a & b\\
c & d
\end{pmatrix}
\tens 
\begin{pmatrix}
a & b\\
c & d
\end{pmatrix},
\quad 
\epsilon 
\begin{pmatrix}
a & b\\
c & d
\end{pmatrix}
=
\begin{pmatrix}
1 & 0\\
0 & 1
\end{pmatrix},
\quad
S 
\begin{pmatrix}
a & b\\
c & d
\end{pmatrix}
=
\begin{pmatrix}
d & -qb\\
-q^{-1}c & a
\end{pmatrix},\]
where we understand $\Delta(a)=a\tens a+b\tens c, \ \epsilon(a)=1,\ S(a)=d,$ etc. By definition, the quantum group $\C_q[SU_2]$ is Hopf algebra $\C_q[SL_2]$ with $q$ real and $\ast$-structure
\[
\begin{pmatrix}
a^* & b^*\\
c^* & d^*
\end{pmatrix}=
\begin{pmatrix}
d & -q^{-1}c\\
-qb & a
\end{pmatrix}.
\]
We used the conventions of \cite{Ma:book} and refer there for the history, which is relevant both to~\cite{Wor} and the Drinfeld theory\cite{Dri}. 

On $\C_q[SU_2],$ there is a connected left-covariant calculus $\Omega^1(\C_q[SU_2])$  in \cite{Wor} with basis, in our conventions,
\[\omega^0=d\extd a-qb\extd c,\quad \omega^+=d\extd b-qb\extd d,\quad \omega^-=qa\extd c- c\extd a\] of left-invariant $1$-forms which is dual to the basis $\{\partial_0,\partial_\pm\}$ of left-invariant vector fields generated by the Chevalley basis $\{H,X_\pm\}$ of $su_2$ (so that $[H,X_\pm]=\pm2X_\pm$ and $[X_+,X_-]=H$). The first order calculus is generated by $\{\omega^0,\omega^\pm\}$ as a left module while the right module structure is given by the bimodule relations
\[\omega^0 f= q^{2|f|}f\omega^0,\quad \omega^{\pm}f=q^{|f|}f\omega^\pm\] for homogeneous $f$ of degree $|f|$ where $|a|=|c|=1,\ |b|=|d|=-1,$ and exterior derivative
\begin{gather*}
\extd a=a\omega^0+q^{-1}b\omega^+,\quad \extd b=-q^{-2}b \omega^0+a\omega^-,\\
\extd c=c\omega^0+q^{-1}d\omega^+,\quad \extd d=-q^{-2} d\omega^0+c\omega^-.
\end{gather*}
These extend to a differential graded algebra $\Omega(\C_q[SU_2])$ that has same dimension as classically. Moreover, it is a $\ast$-DGA with
\[\omega^{0*}=-\omega^0,\quad\omega^{+*}=-q^{-1} \omega^-,\quad \omega^{-*}=-q\omega^+.\]

Since $\C_q[SU_2]$ and $\Omega(\C_q[SU_2])$ are $q$-deformations, from Corollary~\ref{preliecorol}, these must be quantised from some pre-Lie algebra structure of $su_2^*$, which we now compute.  Let $q=e^{\frac{\imath\lambda}{2}}=1+\frac{\imath}{2}\lambda+O(\lambda^2)$ for imaginary number $\lambda.$
The Poisson bracket from the algebra relations is 
\begin{gather*}
\{a,b\}=-\frac{\imath}{2}ab,\quad \{a,c\}=-\frac{\imath}{2}ac,\quad \{a,d\}=-\imath bc,\quad \{b,c\}=0,\\ 
\{b,d\}=-\frac{\imath}{2}bd, \quad\{c,d\}=-\frac{\imath}{2}cd.
\end{gather*}
The reader should not be alarmed with $\imath$, as this is a `complexified' Poisson-bracket on $C^\infty(SU_2,\C)$ and is a real Poisson-bracket on $C^\infty(SU_2,\R)$ when we choose real-valued functions instead of complex-valued functions $a,b,c,d$ here.

As $\extd x=\sum_i(\partial_i x)\omega^i,$ we know
\[\partial_0 \begin{pmatrix}
a & b\\
c & d
\end{pmatrix}=\begin{pmatrix}
a & -b\\
c & -d
\end{pmatrix},\quad \partial_+\begin{pmatrix}
a & b\\
c & d
\end{pmatrix}= \begin{pmatrix}
0 & a\\
0 & c
\end{pmatrix},\quad \partial_-\begin{pmatrix}
a & b\\
c & d
\end{pmatrix}=\begin{pmatrix}
b & 0\\
d & 0
\end{pmatrix}.\]

From $a\omega^0-\omega^0 a=(1-q^2)a\omega^0=-\imath\lambda a\omega^0+O(\lambda^2),$ we know that $\gamma(a,\omega^0)=-\imath a\omega^0.$ Likewise, we can get
\[\gamma(\begin{pmatrix}
a & b\\
c & d
\end{pmatrix},\omega^i)=\frac{1}{2} t_i 
\begin{pmatrix}
a & -b\\
c & -d
\end{pmatrix} \omega^i,\ \forall\,i\in \{0,\pm\},\quad t_0=-2\imath,\ t_\pm=-\imath.\]
Now we can compute the pre-Lie structure $\Xi:su_2^*\tens su_2^*\to su_2^*$ by comparing with (\ref{preconnection-ijk-1}), namely
\[\gamma(\begin{pmatrix}
a & b\\
c & d
\end{pmatrix},\omega^j)=\sum_{i,k\in \{0,\pm\}} \Xi^{ij}_k(\partial_i\begin{pmatrix}
a & b\\
c & d
\end{pmatrix})\omega^k,\] 
we know the only nonzero coefficients are
$\Xi^{00}_0=-\imath,\ \Xi^{0+}_+=-\frac{\imath}{2},\ \Xi^{0-}_-=-\frac{\imath}{2}$. 
Then $\Xi(\phi,\phi)=-\imath \phi,\ \Xi(\phi,\psi_+)=-\frac{\imath}{2}\psi_+,\  \Xi(\phi,\psi_-)=-\frac{\imath}{2}\psi_-,$ and $\Xi$ is zero on other terms, where $\{\phi,\psi_\pm\}$ is the dual basis of $su_2^*$ to $\{H,X_\pm\}$. 
Thus the corresponding pre-Lie structure of $su_2^*$ is 
\begin{gather*}
\Xi (\phi,\phi)=-\imath \phi,\quad \Xi(\phi,\psi_\pm)=-\frac{\imath}{2}\psi_\pm,\quad\text{and zero on all other terms.}
\end{gather*} 
Let $t=-2\imath \phi,\ x_1=\imath (\psi_++\psi_-),\ x_2=\psi_+-\psi_-,$ we have a real pre-Lie structure for $su_2^*=\mathrm{span}\{t,x_1,x_2\}:$
\[t\circ t=-2t,\quad t\circ x_i=-x_i,\ \forall\,i=1,2.\]
This is a $3$D version of $\mathfrak{b}_{1,-2}.$
\end{example}

\begin{example}\label{qt-flat}
Let $\cg$ be a quasi-triangular Lie bialgebra with $r$-matrix $r=r\bo\tens r\bt\in\cg\tens\cg.$ Then $\cg$ acts on its dual $\cg^*$ by coadjoint action $\ad^*$ and by~\cite[Lemma 3.8]{Ma:blie}, $\cg^*$ becomes a left $\cg$-crossed module with $-\Xi,$ where $\Xi$ is left $\cg^*$-action $\Xi(\phi,\psi)=-\<\phi,r\bt\>\ad_{r\bo}\psi.$ To satisfy Poisson-compatibility (\ref{compatible}), $(\cg,r)$ is required to obey $r\bo\tens [r\bt,x]+r\bt\tens [r\bo,x]=0,$ i.e., $r_+\la x=0$ for any $x\in\cg,$ where $r_+=(r+r_{21})/2$ is the symmetric part of $r.$  In this case $\cg^*$ has a pre-Lie algebra structure with $\Xi(\phi,\psi)=-\<\phi,r\bt\>\ad^*_{r\bo}\psi$ by Corollary~\ref{preliecorol}. We see in particular that every finite-dimensional cotriangular Lie bialgebra is canonically a pre-Lie algebra. 
\end{example}

\subsection{Quantum metrics on Example~\ref{2D}} We make a small digression to cover the quantum metrics on the differential calculi in  Example~\ref{2D}.. By definition, given an algebra $A$ equipped with at least $\Omega^1,\Omega^2$ of a DGA, a quantum metric means $g=g\bo\tens g\bt\in \Omega^1\tens_A\Omega^1$ (a formal sum notation) such that $\wedge(g)=0$  (this expresses symmetry via the wedge product) and which also has an inverse $(\ ,\ ):\Omega^1\tens_A\Omega^1\to A$ such that $(\omega,g\bo)g\bt=g\bo(g\bt,\omega)=\omega$ for all $\omega\in\Omega^1$. This data makes $\Omega^1$ left and right self-dual in the monoidal category of $A$-bimodules. It can be shown in this case \cite{BegMa14} that $g$ must be central. In the $*$-DGA case we also require that the metric is `real' in the sense
\[ {\rm flip}(*\tens *)g=g\]
where `flip' interchanges the tensor factors. We analyse the moduli of quantum metrics and also their $\lambda=0$ limit in each of the five cases:

1) For the calculus $\Omega^1(U_\lambda(\mathfrak{b}_{1,\alpha})),$  direct computation shows that there is no nondegenerate metric at the polynomial level. However, if we allow a wider class of functions (including $x^{-1},x^{\pm\alpha}$), then the 1-forms
\[ u=x^{-1}\extd x,\quad v=x^{\alpha}\extd t\]
are central, and obey $u^*=u, v^*=v$ with the result that there is a 3-parameter family of nondegenerate quantum metrics of the form
\begin{equation*}
g=\frac{c_1}{x^2} \extd x\tens \extd x +c_2 x^{\alpha-1}(\extd x\tens\extd t+\extd t\tens \extd x) +c_3 x^{2\alpha}\extd t\tens \extd t,
\end{equation*}
where $c_1,c_2,c_3$ are complex parameters with $\det c=c_1c_3-c_2^2\neq0$ so that $\det g=(\det c) x^{2\alpha-2}\neq 0.$  The metric obeys the `reality' condition if and only if the $c_i$ are real. In the classical limit $\lambda=0$ the metric has constant scalar curvature
\[R=-\frac{2\alpha^2 c_3}{\det c}\]
and with more work one can show that up to a change of coordinates, the metric depends only on the one parameter given by the value of $R$. For example, if we require $\det c<0$ for Minkowski signature, then this metric is essentially de Sitter or anti-de Sitter space for some length scale. This is taken up and extended to the $n$-dimensional case in \cite{MaTao:cos}.

2) For the calculus $\Omega^1(U_\lambda(\mathfrak{b}_{2,\beta\neq0}))$,  one can show that
\[u=x^{\beta-1}\extd x,\quad v=x^{\beta-1}(x\extd t-\beta t\extd x)\] are central $1$-forms. Also $u^*=u,\ v^*=v+\lambda \beta(\beta-2)u.$ Then a quantum metric takes the form
\begin{equation*}
\begin{split}
g=c_1 u\tens u+c_2(u\tens v+v^*\tens u)+c_3(v^*\tens v+\beta\lambda (u\tens v-v^*\tens u))
\end{split}
\end{equation*}
and obeys `reality' condition if an only if the $c_i$ is real. In this case, if we let $t'=t+\frac{c_2}{\beta},$ then $\extd t'=\extd t, v'=x^{\beta-1}(x\extd t'-\beta t'\extd x)=v-c_2u$ and similarly for $v'{}^*$. We can use this freedom to set $c_2=0$ at the expense of a change to $c_1$ and prefer to assume this standard form:
\begin{align*}
g&=c_1 u\tens u+c_3(v^*\tens v+\beta\lambda(u\tens v-v^*\tens u))\\
&=x^{2\beta-2}(c_1+c_3\beta^2t^2-\lambda c_3\beta^2(2\beta-3)t+\lambda^2c_3\beta^2(\beta^2-3\beta+3))\extd x\tens\extd x\\
&-x^{2\beta-1}c_3(t-\lambda \beta(\beta-1))(\extd x\tens \extd t+\extd t\tens \extd x)+x^{2\beta}c_3\extd t\tens \extd t
\end{align*}
where $c_1,c_3\ne 0$. The scaler curvature in the classical limit is 
\[R=-x^{-2\beta} \frac{2 \beta(\beta+1)c_1}{(c_1+c_3(\beta^2-1)t^2)^2}\]
which has a singularity at $x=0$ when $\beta>0$, suggesting some kind of gravitational source. The $\beta=1$ case was already in \cite{BegMa14} and it was shown that for $c_1>0$ and $c_3<0$ the gravitational source at $x=0$ was so strong that even outgoing light was turned back in (or rather any particle of arbitrarily small mass). 

3) For the calculus $\Omega^1(U_\lambda(\mathfrak{b}_{3}))$, there is no quantum metric at the algebraic level, but if we allow a wider class of functions (including $x^{-1}$ and $\ln(x)$), then
\[ u=x^{-1}\extd x,\quad v=x^{-1}\ln (x)\extd x+x^{-1}\extd t\]
are central 1-forms and obey $u^*=u, v^*=v$. Then there is a 3-parameter class of quantum metrics built from these. However, with this wider class of functions, $t'=t+x \ln x$ then takes us back to $\alpha=-1$ in case 1). Here $x,t'$ and their differentials obey the relations in case 1) and $x^{-1}\extd t'=u+v$ here.

4) For the calculus $\Omega^1(U_\lambda(\mathfrak{b}_{4}))$, there is no quantum metric at the polynomial level but if we allow a wider class of functions (including $x^{-1}$), 
\[u=\frac{1}{x^2}\extd t,\quad v=\frac{1}{x^2}(x\extd x-t\extd t)\]
are central $1$-forms  and $u^*=u, v^*=v-3\lambda u.$  A quantum metric then takes the form
\begin{equation*}
g=c_1 u\tens u+c_2(u\tens v+v^*\tens u) + c_3(v^*\tens v+\lambda(u\tens v-v^*\tens u)).
\end{equation*}
Thus $g$ in this form is manifestly central and the `reality' condition requires the $c_i$ to be real. If we choose a new variable $t'=t-c_2$ we can set $c_2=0$ and take our metric in the form 
\begin{align*}
g&=c_1 u\tens u+c_3(v^*\tens v+\lambda(u\tens v-v^*\tens u))\\
&=\frac{c_3}{x^2}\extd x\tens \extd x-c_3\frac{t+2\lambda}{x^3}(\extd x\tens \extd t+\extd t\tens\extd x)+\frac{c_1+c_3(t^2+5\lambda t+\lambda^2)}{x^4}\extd t\tens \extd t,
\end{align*} 
where $c_1,c_3\ne 0$. The Ricci tensor in the classical limit is
\[R=4{x^2-2t^2\over c_1}-{8\over c_3}.\]

5) For  the calculus $\Omega^1(U_\lambda(\mathfrak{b}_{5})),$  we find that 
\[u=\extd x,\quad v=x\extd t+(x-t)\extd x\]
are central $1$-forms and $u^*=u,\ v^*=v-\lambda u.$  Then a general quantum metric takes the form
\begin{equation*}
g=c_1 u\tens u+c_2(u\tens v +v^*\tens u)+c_3(v^*\tens v+\lambda(u\tens v-v^*\tens u)),
\end{equation*} 
which is real when the $c_i$ are real.  In this case if we choose a new variable $t'=t-c_2$ which has the same relations in differential algebra and which can be used to absorb $c_2$ with a different value of $c_3$. We can therefore take our metric in a standard form
\begin{align*}
g&=c_1  u\tens u+c_3(v^*\tens v+\lambda(u\tens v-v^*\tens u))\\
&=(c_1+c_3(t^2-(2x-\lambda)t+x^2+\lambda^2))\extd x\tens\extd x +c_3x(x-t)(\extd x\tens \extd t+\extd t\tens \extd x)\\
&\quad+c_3x^2\extd t\tens \extd t,
\end{align*}
where $c_1,c_3\ne 0$. The Ricci tensor in the classical limit is
\[R=-\frac{4}{c_1x^2}.\]
This case, if we allow a wider class of functions, is equivalent to $\beta=1$ in case 2). Here $t'=t+x\ln x$ obeys the relations there, and $x\extd t'-t'\extd x=x\extd t+(x-t)\extd x$ as here. 

What these various examples show is that the classification of all left-covariant nice differential structures is possible and each has its only highly restrictive form of noncommutative Riemannian geometry associated to it. The quantum-Levi-Civita connection is found for case 1) in \cite{MaTao:cos} and for the $\beta=1$ case of 2) in \cite{BegMa14} and other cases should be achievable in a similar way.

\section{Quantisation of tangent bundle $TG=G\rcross \underline\cg$}

We are interested in quantisation of the tangent bundle $TG$ of a Poisson-Lie group $G$.  The noncommutative coordinate algebra in this case is a bicrossproduct~\cite{Ma:book, Ma:bicross}.

\subsection{Review on bicrossproduct of Hopf algebras}

We start with the notions of double cross sum and bicross-sum of Lie bialgebras~\cite[Ch. 8]{Ma:book}. We say $(\cg,\cm,\ra,\la)$ forms a  right-left matched pair of Lie algebras if $\cg,\cm$ are both Lie algebras and $\cg$ right acts on $\cm$ via $\ra,$ $\cm$ left acts on $\cg$ via $\la$ such that 
\begin{gather*}
[\phi,\psi]\ra\xi=[\phi\ra\xi,\psi]+[\phi,\psi\ra\xi]+\phi\ra(\psi\la\xi)-\psi\ra(\phi\la \xi),\\
\phi\la [\xi,\eta]=[\phi\la\xi,\eta]+[\xi,\phi\la\eta]+(\phi\ra\xi)\la\eta-(\phi\ra\eta)\la\xi,
\end{gather*}
for any $\xi,\eta\in\cg,\, \phi,\psi\in\cm.$ Given such a matched pair, one can define the `double cross sum Lie algebra'  $\cg\dcross\cm$ as the vector space $\cg\oplus\cm$ with the Lie bracket
\[[(\xi,\phi),(\eta,\psi)]=([\xi,\eta]+\phi\la\eta-\psi\la\xi,[\phi,\psi]+\phi\ra\eta-\psi\ra\xi).\]
In addition, if both $\cg$ and $\cm$ are now Lie bialgebras with $\la,\ra$ making $\cg$ a left $\cm$-module Lie coalgebra and $\cm$ a right $\cg$-module Lie coalgebra, such that $\phi\ra \xi\o\tens\xi\t+\phi\o\tens\phi\t\la\xi=0$ for all $\xi\in\cg,\ \phi\in\cm,$ then the direct sum Lie coalgebra structure makes $\cg\dcross\cm$ into a Lie bialgebra, \textit{the double cross sum Lie bialgebra}.

Next, if $\cg$ is finite dimensional, the matched pair of Lie bialgebras $(\cg,\cm,\ra,\la)$ data equivalently defines a \textit{right-left bicross sum Lie bialgebra} $\cm\rlbicross\cg^*$ bulit on $\cm\oplus\cg^*$ with 
\begin{gather}
[(\phi,f),(\psi,h)]=([\phi,\psi]_{\cm},[f,h]_{\cg^*}+f\ra\psi-h\ra\phi),\label{bicross-bra}\\
\delta \phi=\delta_{\cm}\phi+(\id-\tau)\beta(\phi),\quad \delta f=\delta_{\cg^*} f,\label{bicross-cobra}
\end{gather}
for any $\phi,\psi\in\cm$ and $f,h\in\cg^*,$
where the right action of $\cm$ on $\cg^*$ and the left coaction of $\cg^*$ on $\cm$ are induced from $\ra,\la$ by
\[\<f\ra\phi,\xi\>=\<f,\phi\la\xi\>,\quad \beta(\phi)=\sum_i f^i\tens\phi\ra e_i,\]
for all $\phi\in\cm,\,f\in\cg^*,\,\xi\in\cg$ and $\{e_i\}$ is a basis of $\cg$ with dual basis $\{f^i\}.$ We refer to~\cite[Section 8.3]{Ma:book} for the proof. 

Now let $(\cg,\cm,\ra,\la)$ be a matched pair of Lie algebras and $M$ be the connected and simply connected Lie group associated to $\cm.$ The Poisson-Lie group $M\rlbicross \cg^*$ associated to the bicross sum $\cm\rlbicross\cg^*$ is the semidirect product $M\rcross\cg^*$ (where $\cg^*$ is regarded as an abelian group) equipped with Poisson bracket
\[\{f,g\}=0,\quad\{\xi,\eta\}=[\xi,\eta]_\cg,\quad \{\xi,f\}=\alpha_{*\xi}(f)\]
for all $f,g$ functions on $M$ and $\xi,\eta$ linear functions on $\cg^*,$ where $\alpha_{*\xi}$ is the vector field for the action of $\cg$ on $M$. See~\cite[Proposition 8.4.7]{Ma:book} for the proof.
Note that here $\cg,\cm$ are both viewed as Lie bialgebras with zero cobracket, so the Lie bracket and Lie cobracket of the bicross sum Lie bialgebra $\cm\rlbicross\cg^*$ here is given by (\ref{bicross-bra}) and (\ref{bicross-cobra})  but with $[\ ,\ ]_{\cg^*}=0,\delta_\cm=0$ there.

More precisely, let $(\cg,\cm,\ra,\la)$ be a matched pair of Lie algebras, with the associated connected and simply connected Lie groups $G$ acting on $\cm$ and $M$ acting on $\cg.$ The action $\ra$ can be viewed as Lie algebra cocycle $\ra\in Z^1_{\la^*\tens {-\id}}(\cm,\cg^*\tens\cm)$ and under some assumptions then can be exponentiated to a group cocycle $a\in Z^1_{\la^*\tens \Ad_R}(M,\cg^*\tens\cm),$
which defines an infinitesimal action of $\cg$ on $M$. Hence, by evaluation of the corresponding vector fields, $a$ defines a left action of the Lie algebra $\cg$ on $C^\infty(M)$~\cite{Ma:mat}:
\begin{equation}\label{vfieldaction}
(\widetilde{\xi} f)(s)=\widetilde{a_\xi} (f)(s)=\frac{\extd}{\extd t}\Big|_{{t=0}} f(s \exp(t a_\xi(s))),\quad\forall\,f\in C^\infty(M),\ \xi\in\cg.
\end{equation} 
We also note that $\cm$ acts on $M$ by left invariant vector field $(\widetilde{\phi}f)(s)=\frac{\extd}{\extd t}\Big|_{t=0} f(s\exp{(t\phi)})$ for any $\phi\in \cm,\,f\in C^\infty(M)$ and these two actions fit together to an action of $\cg\dcross\cm$ on $C^\infty(M).$ 

Finally,  we can explain the bicrossproduct $\C[M]\lrbicross U_\lambda(\cg)$ based on a matched pair of Lie algebras $(\cg,\cm,\ra,\la)$, where $\C[M]$ is the algebraic model of functions on $M.$ The algebra of $\C[M]\lrbicross U_\lambda(\cg)$ is the cross product defined by the action (\ref{vfieldaction}). Its coalgebra, on the other hand, is the cross coproduct given in reasonable cases by the right coaction (defined by the left action of $M$ on $\cg$)
\[\beta:\cg\to\cg\tens \C[M],\quad \beta(\xi)(s)=s\la\xi,\quad\ \forall\,\xi\in\cg,\ s\in M.\] The map $\beta$ is extended to products of the generators of $U_\lambda(\cg)$ in aim to form a bicrossproduct $\C[M]\lrbicross U_\lambda(\cg)$ as in \cite[Theorem 6.2.2]{Ma:book}. 

The Poisson-Lie group $M\rlbicross\cg^*$ quantises to  $\C[M]\lrbicross U_\lambda(\cg)$ as a noncommutative deformation of the commutative algebra of functions $\C[M\rlbicross \cg^*]\subset C^\infty(M\rlbicross\cg^*)$.  See more details in~\cite[Section 8.3]{Ma:book}. The half-dualisation process we have described at the Lie bialgebra level also works at the Hopf algebra level, at least in the  finite-dimensional case. So morally speaking, $U_\lambda(\cg)\dcross U(\cm)$ half-dualises in a similar way to the bicrossproduct Hopf algebras $\C[M]\lrbicross U_\lambda(\cg).$ If one is only interested in the algebra and its calculus, we can extend to cross product $C^\infty(M)\lcross U_\lambda(\cg).$

\subsection{Poisson-Lie group structures on tangent bundle $G\rcross\underline\cg$} Let $G$ be a Lie group with Lie algebra $\cg.$ As a Lie group, the tangent bundle $TG$ of Lie group $G$ can be identified with the semidirect product of Lie groups $G\rcross\underline\cg$ (under right adjoint action of $G$ on $\cg$) with product 
\[(g_{_{1}},x)(g_{_{2}},y)=(g_{_{1}}g_{_{2}},\Ad(g_{_{2}}^{-1})(x)+y),\quad\forall\,g_{_{1}},g_{_{2}}\in\,G,\ x,y\in\,\cg,\]
where $\underline{\cg}$ is $\cg$ but viewed as an abelian Poisson-Lie group under addition. Naturally, the Lie algebra of $G\rcross\underline\cg$ is the semidirect sum Lie algebra $\cg\rcross\underline{\cg}$ with Lie bracket
\[[\xi,\eta]=[\xi,\eta]_{\cg},\quad [x,y]=0,\quad [x,\xi]=[x,\xi]_\cg,\quad\forall\,\xi,\eta\in\cg,\,x,y\in\underline{\cg}.\]

In light of the observations in Section~5.1, we propose the following construction on the Poisson-Lie structure on tangent bundle $G\rcross\underline{\cg}$ via bicross sum. In what follows we assume that $G$ is a finite dimensional connected and simply connected Poisson-Lie group, and $\cg$ is its Lie algebra with the corresponding Lie bialgebra structure.

Denote by $\overline{\cg^*}:=(\cg^*,[\ , \ ]_{\cg^*},\text{ zero Lie cobracket})$  and $\overline{\cg}:=(\cg, [\ ,\ ]_\cg,\text{ zero Lie cobracket}),$ where $\overline{\cg^*}$ is the dual of  Lie bialgebra $\underline{\cg}=(\cg,\text{ zero bracket},\delta_\cg)$. One can check that $\overline{\cg^*}$ and $\overline{\cg}$ form together a matched pair of Lie bialgebra with coadjoint actions, i.e., 
\[\xi\ra \phi=-\ad^*_\phi\xi=\<\phi,\xi\o\>\xi\t,\quad \xi\la \phi=\ad^*_\xi \phi=\phi\o\<\phi\t,\xi\>\]
for any $\phi\in \overline{\cg^*},\,\xi\in\overline{\cg}.$ 

\subsubsection{Lie bialgebra level} The double cross sum Lie bialgebra $\overline{\cg^*}\dcross\overline{\cg}$ is then built on $\cg^*\oplus\cg$ as a vector space with Lie bracket
\begin{gather*}
[\phi,\psi]=[\phi,\psi]_{\cg^*},\quad [\xi,\eta]=[\xi,\eta]_{\cg},\quad \forall\,\phi,\psi\in\overline{\cg^*},\,\xi,\eta\in\overline{\cg},\\
[\xi,\phi]=\xi\ra\phi+\xi\la\phi=\<\phi,\xi\o\>\xi\t+\phi\o\<\phi\t,\xi\>,
\end{gather*}
and zero Lie cobracket. This is nothing but the Lie algebra of $D(\cg)=\cg^*\dcross\cg$ (the double of Lie bialgebra $\cg$) with zero Lie-cobracket.

Correspondingly, the right-left bicross sum Lie bialgebra defined by the matched pair $(\overline{\cg^*},\overline{\cg},\ra,\la)$ above is $\overline{\cg}\rlbicross\underline{\cg}$, whose Lie algebra is semidirect sum $\overline{\cg}\rcross\underline{\cg}$ and the Lie coalgebra is semidirect cobracket $\overline{\cg}\lcocross\underline{\cg},$ namely
\begin{equation}\label{LB-T}
\begin{split}
[\xi,\eta]=[\xi,\eta]_{\cg},\quad [x,y]=0,\quad [x,\xi]=[x,\xi]_\cg;\quad\textbf{     }\\
\delta\xi=(\id-\tau)\delta_\cg(\xi)=\underline{\xi\o}\tens\overline{\xi\t}-\overline{\xi\t}\tens \underline{\xi\o},\quad \delta x=\delta_\cg x,\quad
\end{split}
\end{equation}
for any $\xi,\eta\in\overline{\cg},\,x,y\in\underline{\cg}$. Here the coaction on $\overline{\cg}$ is the Lie cobracket $\delta_\cg$ viewed as map from $\overline{\cg}$ to $\underline{\cg}\tens\overline{\cg}.$

\subsubsection{Poisson-Lie level} Associated to the right-left bicross sum Lie bialgebra $\overline{\cg}\rlbicross\underline{\cg},$ the Lie group $G\rcross \underline{\cg}$ is a Poisson-Lie group (denoted by $\overline{G}\rlbicross\underline{\cg}$) with the Poisson bracket
\begin{equation}\label{PB-T}
\{f,h\}=0,\quad \{\phi,\psi\}=[\phi,\psi]_{\cg^*},\quad \{\phi,f\}=\widetilde{\phi}f,
\end{equation}
for any $\phi,\psi\in \overline{\cg^*}\subseteq C^\infty(\underline{\cg})$ and $f,h\in C^\infty(\overline{G}),$ 
where $\widetilde{\phi}$ denotes the left Lie algebra action of $\overline{\cg^*}$ on $C^\infty(G)$ (viewed as vector field on $G$) and is defined by the right action of $\overline{\cg^*}$ on $\overline{\cg}.$ 

The vector field $\widetilde{\phi}$ for any $\phi\in\cg^*$ in this case can be interpreted more precisely. We can view the actions between $\overline{\cg^*}$ and $\overline{\cg}$ as Lie algebra $1$-cocycles, namely the right coadjoint action $\ra=-\ad^*:\overline{\cg}\tens \overline{\cg^*}\to \overline{\cg}$ (of $\overline{\cg^*}$ on $\overline{\cg}$) is viewed as a map $\overline{\cg}\to (\overline{\cg^*})^*\tens\overline{\cg}= (\underline{\cg})^{**}\tens\overline{\cg}=\underline{\cg}\tens\overline{\cg}.$ It maps $\xi$ to $\sum_i e_i\tens\xi\ra f^i=\sum_i e_i\tens \<f^i,\xi\o\>\xi\t=\xi\o\tens\xi\t,$ which is nothing but the Lie cobracket $\delta_\cg$ of $\cg.$ Likewise, the left coadjoint action of $\overline{\cg}$ on $\overline{\cg^*}$ is viewed as the Lie cobracket $\delta_{\cg^*}$ of $\cg^*.$ We already know that Lie $1$-cocycle $\delta_\cg\in Z^1_{-\ad}(\cg,\cg\tens\cg)$  exponentiates to group cocycle
\[D^\vee\in Z^1_{\Ad_R}(G,\cg\tens\cg),\]
thus
\begin{equation}\label{phitilde}
\widetilde{\phi}_g:=(L_{g})_*((\phi\tens\id)D^\vee(g))\in T_gG,\quad\forall\,g\in G,
\end{equation}
defines the vector field on $G.$ 

According to \cite[Proposition 8.4.7]{Ma:book},  the Poisson bivector on tangent bundle $TG=\overline{G}\rlbicross\underline{\cg}$ is
\begin{equation}\label{PV-T}
P=\sum_i(\partial_i\tens \widetilde{f^i}-\widetilde{f^i}\tens\partial_i)+\sum_{i,j,k} d^{ij}_k f^k \partial_i\tens \partial_j
\end{equation}
where $\{\partial_i\}$ is the basis of left-invariant vector fields generated by the basis $\{e_i\}$ of $\cg$ and $\{f^i\}$ is the dual basis of $\cg^*.$ Here  $P_{KK}=\sum_{i,j,k} d^{ij}_k f^k \partial_i\tens \partial_j$ is the known Kirillov-Kostant bracket on $\underline\cg$ with $\delta_\cg e_k=\sum_{ij}d^{ij}_k e_i\tens e_j.$ 
We arrive at the following special case of \cite[Proposition~8.4.7]{Ma:book}.

\begin{lemma}\label{tangent}
Let $G$ be a finite dimensional connected and simply connected Poisson-Lie group and $\cg$ be its Lie algebra.
The tangent bundle $TG={G}\rcross\underline{\cg}$ of $G$ admits a Poisson-Lie structure given by (\ref{PB-T}) or (\ref{PV-T}), denoted by $\overline{G}\rlbicross\underline{\cg}$. The corresponding Lie bialgebra is $\overline{\cg}\rlbicross\underline{\cg}$ given by (\ref{LB-T}).
\end{lemma}

\subsubsection{Bicrossproduct Hopf algebra} Finally,  when the actions and coactions are suitably algebraic, we have a bicrossproduct Hopf algebra $\C[\overline{G}]\lrbicross U_\lambda(\overline{\cg^*})$ as a quantisation of the commutative algebra of functions $\C[\overline{G}\rlbicross \underline{\cg}]$ on tangent bundle $\overline{G}\rlbicross \underline{\cg}$ of Poisson-Lie group $G.$ The commutation relations of $\C[\overline{G}]\lrbicross U_\lambda(\overline{\cg^*})$ are
\[[f,h]=0,\quad [\phi,\psi]=\lambda[\phi,\psi]_{\cg^*},\quad [\phi,f]=\lambda\widetilde{\phi} f,\]
for any $\phi,\psi\in \overline{\cg^*}\subseteq C^\infty(\underline{\cg})$ and $f,h\in \C[\overline{G}].$ This construction is still quite general but includes a canonical example for all compact real forms $\cg$ of complex simple Lie algebras based in the Iwasawa decomposition to provide the double cross product or `Manin triple' in this case\cite{Ma:mat}. We start with an even simpler example.

\begin{example}\label{Tm^*}
Let $\cm$ be a finite dimensional real Lie algebra, viewed as a Lie bialgebra with zero Lie-cobracket. Take $G=\cm^*$ the abelian Poisson-Lie group with Kirillov-Kostant Poisson bracket given by the Lie bracket of $\cm.$ Then $\cg=\cm^*$ and $\overline{\cg^*}=\cm$ and $\overline{\cg}=\overline{\cm^*}=\mathbb{R}^n,$ where $n=\dim\cm.$ Since the Lie bracket of $\cm^*$ is zero, $\overline{\cm^*}$ acts trivially on $\cm$, while $\cm$ acts on $\overline{\cm^*}$ by right coadjoint action $-\ad^*$, namely, 
\[f\ra\xi=-\ad^*_\xi f,\quad\text{or}\ \<f\ra\xi,\eta\>=\<f,[\xi,\eta]_{\cm}\>\] 
for any $f\in\overline{\cm^*},\,\xi,\eta\in\cm$. So $(\cm,\overline{\cm^*},\ra=-\ad^*,\la=0)$ forms a matched pair.  

The double cross sum of the matched pair $(\cm,\overline{\cm^*})$ is $\cm\rcross\overline{\cm^*},$ the semidirect sum Lie algebra with coadjoint action of $\cm$ on $\overline{\cm^*}$
\begin{gather*}
[\xi,\eta]=[\xi,\eta]_\cm,\quad [f,h]=0,\quad [f,\xi]=f\ra\xi=\<\xi,f\o\>f\t,\\
\delta\xi=0,\quad\delta f=0,\quad\forall\,\xi,\eta\in\cm,\,f,h\in\overline{\cm^*}.
\end{gather*}
Meanwhile, the right-left bicross sum of the matched pair $(\cm,\overline{\cm^*})$ is $\overline{\cm^*}\lcocross\cm^*,$ the semidirect sum Lie coalgebra
\begin{gather*}
[f,h]=0,\quad [\phi,\psi]=0,\quad [\phi,f]=\phi\ra f=0,\\
\delta f=(\id-\tau)\beta(f),\quad \delta\phi=\delta_{\cm^*}\phi,
\end{gather*}
for any $f,h\in \overline{\cm^*},\,\phi,\psi\in\cm^*,$ where the left coaction of $\cm^*$ on $\overline{\cm^*}$ is given by
\[\beta:\overline{\cm^*}\to\cm^*\tens\overline{\cm^*},\quad \beta(f)=\sum_i f^i\tens f\ra e_i,\] and $\{e_i\}$ is a basis of $\cm$ with dual basis $\{f^i\}$ of $\cm^*.$

The tangent bundle of $\cm^*$ is the associated Poisson-Lie group of $\overline{\cm^*}\lcocross\cm^*$, which is $\overline{M^*}\lcocross\cm^*=\R^n\lcocross\cm^*,$ an abelian Lie group, where we identify the abelian Lie group $\overline{M^*}$ with its abelian Lie algebra $\overline{\cm^*}.$  Let $\{x^i\}$ be the coordinate functions on $\R^n$ identified with $\{e_i\}\subset\cm\subseteq C^\infty(\overline{\cm^*})=C^\infty(\R^n),$ as $e_i(\sum_j\lambda_jf^j)=\lambda_i.$ 
The right action of $\cm$ on $\overline{\cm^*}$ transfers to $\delta_{\cm^*}\in Z^1(\overline{\cm^*},\underline{\cm^*}\tens \overline{\cm^*}).$ As Lie group $\overline{M^*}$ is abelian and $\overline{M^*}=\overline{\cm^*}=\R^n,$ so the associated group cocycle is identical to $\delta_{\cm^*},$ thus from (\ref{phitilde}), we have
$\widetilde{\xi}_x f =\<x\o,\xi\>x\t{}_x f=\sum_i\<x\o,\xi\>\<x\t,e_i\>f^i{}_x f=\sum_i\<[\xi,e_i]_\cm,x\>\frac{\partial f}{\partial x^i}(x),$
where we use the Lie cobracket in an explicit notation. This shows that \[\widetilde{\xi}=\sum_{i,j,k}\<f^i,\xi\> c^k_{ij}x^k\frac{\partial}{\partial x^j},\quad\forall\,\xi\in\cm,\]
where $c^k_{ij}$ are the structure coefficients of Lie algebra $\cm,$ i.e., $[e_i,e_j]_\cm=\sum_k c^k_{ij}e_k.$ Therefore the Poisson-bracket on $\R^n\lcocross\cm^*$ is given by
\[\{f,h\}=0,\quad\{\xi,\eta\}=[\xi,\eta]_\cm,\quad \{\xi,f\}=\widetilde{\xi} f=\sum_{i,j,k}\<f^i,\xi\> c^k_{ij}x^k\frac{\partial f}{\partial x^j},\] where $f,h\in C^\infty(\R^n),\,\xi,\eta\in\cm.$

The bicrossproduct Hopf algebra $\C[\overline{G}]\lrbicross U_\lambda(\overline{\cg^*})=\C[\R^n]\lcross U_\lambda(\cm),$ as the quatisation of $C^\infty(\R^n\lcocross\cm^*),$ has commutation relations
\[[x^i,x^j]=0,\quad[e_i,e_j]=\lambda\sum_{k} c^k_{ij}e_k,\quad [e_i,x^j]=\lambda \sum_k c^k_{ij}x^k\]
where $\{x^i\}$ are coordinate functions of $\R^n=\overline{\cm^*}$, identified with $\{e_i\}$ basis of $\cm.$ As an algebra we can equally well take $C^\infty(\R^n)\lcross U_\lambda(\cm)$, i.e., not limiting ourselves to polynomials. Then $[e_i,f]=\lambda \sum_{j,k} c^k_{ij}x^k\frac{\partial f}{\partial x^j}$ more generally for the cross relations. 
\end{example}

\begin{example}\label{TSU_2} We take $SU_2$ with its standard Lie bialgebra structure on $su_2$, where the matched pair comes from the Iwasawa decomposition of $SL_2(\C)$\cite{Ma:mat}. 
 The bicrossproduct Hopf algebra $\C[{SU_2}]\lrbicross U_\lambda({su^*_2}),$ as an algebra, is the cross product $\C[SU_2]\lcross U_\lambda(su_2^*)$ with $a,b,c,d$ commutes, $ad-bc=1$, $[x^i,x^3]=\lambda x^i$ ($i=1,2$) and  
\[[x^i,\mathbf{t}]=\lambda\mathbf{t}[e_i,\mathbf{t}^{-1}e_3\mathbf{t}-e_3],\ i=1,2,3,\] namely,
\begin{gather}
[x^1,\mathbf{t}]=-\lambda bc\,\mathbf{t}e_2+\frac{\lambda}{2}\mathbf{t}\,\mathrm{diag}(ac,-bd)+\frac{\lambda}{2}\mathrm{diag}(b,-c),\nonumber\\
[x^2,\mathbf{t}]=\lambda bc\,\mathbf{t} e_1-\frac{\imath \lambda}{2}\mathbf{t}\,\mathrm{diag}(ac,bd)+\frac{\imath \lambda}{2}\,\mathrm{diag}(b,c),\label{TSU_2-action}\\
[x^3,\mathbf{t}]=-\lambda ad\,\mathbf{t}+\lambda\mathrm{diag}(a,d),\nonumber
\end{gather}
where we denote $\mathbf{t}=\begin{pmatrix} a & b\\ c & d \end{pmatrix}$ and $\{e_i\}$ and $\{x^i\}$ are bases of $su_2$ and $su_2^*$ as the half-real forms of $sl_2(\C)$ and $sl_2^*(\C)$ respectively.
The coalgebra of $\C[{SU_2}]\lrbicross U_\lambda({su^*_2})$
is the cross coproduct $\C[SU_2]\rcocross U_\lambda(su_2^*)$ associated with
\begin{gather*}
\Delta(x^i)=1\tens x^i-2\sum_k x^k\tens \mathrm{Tr}(\mathbf{t} e_i \mathbf{t}^{-1} e_k),\quad \epsilon (x^i)=0,\quad\forall\,i\in\{1,2,3\}.
\end{gather*}
The $\ast$-structure is the known one on $\C[SU_2]$ with $x^{i*}=-x^i$ for each $i.$
\end{example}
\proof We already know that the coordinate algebra $\C[SU_2]$ is the commutative algebra $\C[a,b,c,d]$ modulo the relation $ad-bc=1$ with $\ast$-structure $\begin{pmatrix} a^* & b^*\\ c^* & d^* \end{pmatrix} =\begin{pmatrix} d & -c\\ -b & a \end{pmatrix}.$ As a Hopf $\ast$-algebra, the coproduct, counit and antipode of $\C[SU_2]$ are given by
$\Delta \begin{pmatrix} a & b\\ c & d \end{pmatrix} =
\begin{pmatrix} a & b\\ c & d \end{pmatrix} \tens 
\begin{pmatrix} a & b\\ c & d \end{pmatrix},
\ 
\epsilon \begin{pmatrix} a & b\\ c & d \end{pmatrix}=
\begin{pmatrix} 1 & 0\\ 0 & 1 \end{pmatrix},
\ 
S \begin{pmatrix}a & b\\ c & d \end{pmatrix}
= \begin{pmatrix} d & -b\\ -c & a \end{pmatrix}.$

Let $\{H,X_\pm\}$ and $\{\phi,\psi_\pm\}$ be the dual bases of $sl_2(\C)$ and $sl_2^*(\C)$ respectively, where $H=\begin{pmatrix}1 & 0\\ 0 & -1 \end{pmatrix}
,$ $X_+=\begin{pmatrix}0 & 1\\ 0 & 0 \end{pmatrix}
$ and $X_-=\begin{pmatrix}0 & 0\\ 1 & 0 \end{pmatrix}
.$ As the half-real forms of $sl_2(\C)$ and $sl^*_2(\C)$, the Lie algebras $su_2$ and $su^*_2$ have bases
\begin{gather*} 
e_1=-\frac{\imath}{2}(X_++X_-),\quad e_2=-\frac{1}{2}(X_+-X_-),\quad e_3=-\frac{\imath}{2}H,\\
x^1=\psi_++\psi_-,\quad x^2=\imath(\psi_+-\psi_-),\quad x^3=2\phi.
\end{gather*} respectively.
Note that $x^i=-\imath f^i$ where $\{f^i\}$ is the dual basis of $\{e_i\}.$

The Lie brackets and Lie cobrackets of $su_2$ and $su^*_2$ are given by
\begin{gather*}
[e_i,e_j]=\epsilon_{ijk}e_k,\quad \delta e_i=\imath e_i\wedge e_3,\ \forall\,i,j,k;\\
[x^1,x^2]=0,\quad [x^i,x^3]=x^i,\ i=1,2,\quad
\delta x^1=\imath (x^2\tens x^3-x^3\tens x^2),\\
\delta x^2=\imath (x^3\tens x^1-x^1\tens x^3),\quad
\delta x^3=\imath (x^1\tens x^2-x^2\tens x^1),
\end{gather*} where $\epsilon_{ijk}$ is totally antisymmetric and $\epsilon_{123}=1.$
Writing $\xi=\xi^ie_i\in su_2$ and $\phi=\phi_ix^i\in su^*_2$ for $3$-vectors $\vec{\xi}=(\xi^i),\,\vec{\phi}=(\phi_i),$ we know that $(su^*_2,su_2)$ forms a the matched pair of Lie bialgebra with interacting actions 
\[\vec{\xi}\ra \vec{\phi}=(\vec{\xi}\times\vec{e_3})\times\vec{\phi},\quad\vec{\xi}\la\vec{\phi}=\vec{\xi}\times \vec{\phi}.\]

To obtain the action of $su_2^*$ on $\C[SU_2],$ we need to solve 
\cite[Proposition 8.3.14]{Ma:book}
\[\frac{\extd}{\extd t}\Big|_0 a_\phi(\mathrm{e}^{t\xi}u)=\Ad_{u^{-1}}(\xi\ra(u\la\phi)),\quad a_\phi(I_2)=0.\] Note that $SU_2$ acts on $su_2^*$ by $u\la\vec{\phi}=\mathrm{Rot}_u\vec{\phi},$ where we view $\phi$ as an element in $su_2$ via $\rho(\phi)=\phi_ie_i.$ One can check that \[a_{\vec{\phi}}(u)=\vec{\phi}\times\left(\mathrm{Rot}_{u^{-1}}(\vec{e_3})-\vec{e_3}\right)\] is the unique solution to the differential equation. Now we can compute by (\ref{vfieldaction})
\begin{align*}
(\phi\la \mathbf{t}^i{}_j)(u)&=\frac{\extd}{\extd t}\Big|_0 \mathbf{t}^i{}_j (u \mathrm{e}^{t a_\phi(u)})\\
&=\sum_k\frac{\extd}{\extd t}\Big|_0 \mathbf{t}^i{}_k(u) \mathbf{t}^k{}_j(\mathrm{e}^{t a_\phi(u)})\\
&=\sum_k u^i{}_k (a_\phi(u))^k{}_j\\
&=\sum_k u^i{}_k [\rho(\phi),u^{-1}e_3 u-e_3]^k{}_j,
\end{align*}
where $\rho(\phi)=\sum_i\phi_ie_i.$ 
This shows that 
\[[x^i,\mathbf{t}]=\lambda x^i\la\mathbf{t}=\lambda\mathbf{t}[e_i,\mathbf{t}^{-1}e_3\mathbf{t}-e_3]\] 
as displayed. For each $i,$ we can work out the terms on  the right explicitly (using $ad-bc=1$) as
\begin{gather*}
[x^1,\mathbf{t}]=-\frac{\lambda}{2}\begin{pmatrix} abd-a^2c-2b, & b^2d-a^2d+a\\ ad^2-ac^2-d, & bd^2-acd+2c \end{pmatrix},\\
[x^2,\mathbf{t}]=-\frac{\imath \lambda}{2}\begin{pmatrix} a^2c+abd-2b, & a^2d+b^2d-a\\ ac^2+ad^2-d, & bd^2+acd-2c \end{pmatrix},\\
[x^3,\mathbf{t}]=-\lambda \begin{pmatrix} a^2d-a, & abd\\ acd, & ad^2-d \end{pmatrix}.
\end{gather*}
These can be rewritten as the formulae (\ref{TSU_2-action}) we stated.

For convenience, we use Pauli matrices $\sigma_1=\begin{pmatrix} 0 & 1\\ 1 & 0 \end{pmatrix},\ \sigma_2=\begin{pmatrix} 0 & -\imath\\ \imath & 0 \end{pmatrix},\ \sigma_3=\begin{pmatrix} 1 & 0\\ 0 & -1 \end{pmatrix}.$ Clearly, $e_i=-\frac{\imath}{2}\sigma_i$ and $\sigma_i$ obey identities like $\sigma_i\sigma_j=\delta_{ij}I_2+\imath \epsilon_{ijk}\sigma_k$ and $[\sigma_i,\sigma_j]=2\imath \epsilon_{ijk}\sigma_k.$ 

The coaction of $\C[SU_2]$ on $su^*_2$ is defined by $\beta(\phi)(u)=u\la\phi=\mathrm{Rot}_u\vec{\phi}$ for any $u\in SU_2,\,\phi\in su^*_2.$ Again, we view $\phi$ as an element in $su_2,$ so $\rho(u\la\phi)=u\rho(\phi)u^{-1},$ namely $\sum_i(u\la\phi)_i\sigma_i=\sum_i \phi_i u\sigma_i u^{-1}.$ In particular, we have
\[(u\la x^i)_1\sigma_1+(u\la x^i)_2\sigma_2+(u\la x^i)_3\sigma_3=u\sigma_i u^{-1},\quad i=1,2,3.\] 
Note that $\sigma_i\sigma_j=\delta_{ij}I_2+\imath\epsilon_{ijk}\sigma_k$ and $\mathrm{Tr}(\sigma_i\sigma_j)=2\delta_{ij}.$
Time $\sigma_k$ from right to the displayed equation above and then take trace on both sides, we have $2 (u\la x^i)_k=\mathrm{Tr}(u\sigma_i u^{-1}\sigma_k).$ Therefore $u\la x^i=\frac{1}{2}\sum_k\mathrm{Tr}(u\sigma_iu^{-1}\sigma_k)x^k=-2\sum_k\mathrm{Tr}(ue_iu^{-1}e_k)x^k$ and thus $\beta(x^i)=\frac{1}{2}\sum_k x^k\tens \mathrm{Tr}(\mathbf{t}\sigma_i \mathbf{t}^{-1}\sigma_k)=-2\sum_k x^k\tens \mathrm{Tr}(\mathbf{t}e_i \mathbf{t}^{-1}e_k)$. This gives rise to the coproduct of $x^i$ as stated.
\endproof

\subsection{Preconnections on the tangent bundle $\overline{G}\rlbicross \underline{g}$. }

We use the following lemma to construct left pre-Lie structures structure of  $(\overline{\cg}\rlbicross\underline{\cg})^*=(\overline{\cg})^*\lrbicross (\underline{\cg})^*=\underline{\cg^*}\lrbicross \overline{\cg^*},$ whose Lie bracket is the semidirect sum $\underline{\cg^*}\lcross\cg^*$ and the Lie cobracket is the semidirect cobracket ${\cg^*}\rcocross\overline{\cg^*},$ namely
\begin{gather*}
[\phi,\psi]=0,\quad [f,\phi]=f\la\phi=[f,\phi]_{\cg^*},\quad [f,g]=[f,g]_{\cg^*},\\
\delta\phi=\delta_{\cg^*}\phi=\phi\o\tens\phi\t,\quad\delta f=(f\o,0)\tens (0,f\t)-(f\t,0)\tens(0,f\o),
\end{gather*} for any $\phi,\psi\in\underline{\cg^*},\ f,g\in\overline{\cg^*}.$ For convenience, we denote $f\in\cg^*$ by $\overline{f}$ (or $\underline{f}$ ) if viewed in $\overline{\cg^*}$ (or $\underline{\cg^*}$). Thus $\delta f=\underline{f\o}\tens\overline{f\t}+\overline{f\o}\tens\underline{f\t}$ for any $f\in\overline{\cg^*}.$

\begin{lemma}\label{semipre}
Let $(A,\circ)$ be a left pre-Lie algebra and $(B,\ast)$ be a left pre-Lie algebra in $\cg_{_{A}}\mathcal{M},$ namely, there is a left $\cg_{_{A}}$-action $\la$ on $B$ such that
\begin{equation}\label{semipre-con}
a\la (x\ast y)=(a\la x)\ast y+ x\ast (a\la y),
\end{equation}
for any $a,b\in A,\,x,y\in B.$
Then there is  a left pre-Lie algebra structure on $B\oplus A$
\[(x,a)\tilde{\circ} (y,b)=(x\ast y+a\la y,a\circ b).\] Denote this pre-Lie algebra by $B\lcross A$, we have $\cg_{_{B\rtimes A}}=\cg_{_{B}}\lcross \cg_{_{A}}.$
\end{lemma}
\proof Checking the definition of left pre-Lie algebra directly. \endproof

\begin{corollary}\label{semipre-m}
Let $(\cm,\circ)$ be a left pre-Lie algebra. Suppose it admits a (not necessarily unital) commutative associative product $\cdot$ such that
\[[\xi,x\cdot y]_\cm=[\xi,x]_\cm\cdot y+x\cdot [\xi,\eta]_\cm,\quad \forall\,\xi,x,y\in\cm,\]
where $[\ ,\ ]_\cm$ is the Lie bracket defined by $\circ.$ Denote the underlying pre-Lie algebra $\underline{\cm}=(\cm,\cdot).$
Then $\underline{\cm}\lcross_{\ad}\cm$ is a left pre-Lie algebra with product
\begin{equation}
(x,\xi)\tilde{\circ} (y,\eta)=(x\cdot y+[\xi,y]_\cm, \xi\circ \eta)
\end{equation} 
for any $x,y\in\underline{\cm},\,\xi,\eta\in\cm.$
\end{corollary}
\proof Take $(A,\circ)=(\cm,\circ)$ and $(B,\ast)=(\cm,\cdot)$ in Lemma~\ref{semipre}. Here $(\cm,\circ)$ left acts on $(\cm,\cdot)$ by adjoint action and (\ref{semipre-con}) is exactly the condition displayed. \endproof

The assumption made in Corollary~\ref{semipre-m} is that $(\cm,\cdot, [\ ,\ ])$ is a (not necessarily unital) Poisson algebra with respect to the Lie bracket, and that the latter admits a left pre-Lie structure $\circ.$  
 
\begin{theorem}\label{prelie-T} Let $G$ be a finite-dimensional connected and simply connected Poisson-Lie group with Lie bialgebra $\cg.$ Assume that $(\cg^*,[\ ,\ ]_{\cg^*})$ admits a pre-Lie structure $\circ$ and also that $\cg^*$ admits a (not necessarily unital) Poisson algebra structure $(\cg^*,\ast, [\ ,\ ]_{\cg^*})$ 
\begin{equation}\label{pa}
[f,\phi\ast\psi]_{\cg^*}=[f,\phi]_{\cg^*}\ast\psi+\phi\ast [f,\psi]_{\cg^*},
\end{equation}
for any $\phi,\psi\in\underline{\cg^*},\,f\in \cg^*.$ 
Then the semidirect sum $\underline{\cg^*}\lcross \cg^*$ admits a pre-Lie algebra product $\tilde{\circ}$
\begin{equation}\label{circtilde-T}
(\phi,f)\tilde{\circ}(\psi,h)=(\phi\ast\psi+[f,\psi]_{\cg^*},f\circ h)
\end{equation}
and the tangent bundle $\overline{G}\rlbicross\underline{\cg}$ in Lemma~\ref{tangent} admits a Poisson-compatible left-covariant flat preconnection.
\end{theorem}
\proof Take $\cm=\cg^*$ in Corollary~\ref{semipre-m}. We know $\underline{\cg^*}\lcross\cg^*$ is the Lie algebra of $\underline{\cg^*}\lrbicross \overline{\cg^*},$ dual to Lie algebra $\overline{\cg}\rlbicross \underline{\cg}$ of tangent bundle. 
Then we apply Corollary~\ref{preliecorol}.  
\endproof

The corresponding preconnection can be computed by (\ref{preconnection-ijk-1}) explicitly. 
For Poisson-Lie group $G,$ let $\{e_i\}$ be a basis of $\cg$ and $\{f^i\}$ be dual basis of $\cg^*.$ Denote $\{\omega^i\}$ the basis of left--invariant $1$-forms that is dual to $\{\partial_i\}$ the left-invariant vector fields of $G$ generated by $\{e_i\}$ as before. 
For the abelian Poisson-Lie group $\underline\cg$ with Kirillov-Kostant Poisson bracket, let $\{E_i\}$ be a basis of $\underline{\cg}$ and $\{x^i\}$ be the dual basis of $\overline{\cg^*}.$ 
Then $\{ \extd x^i\}$ is the basis of left-invariant $1$-forms that is dual to $\{ \frac{\partial}{\partial x^i}\}$ the basis of the left-invariant vector fields of $\underline\cg$, which are generated by 
$\{E_i\}.$ Now we can choose $\{e_i, E_i\}$ to be the basis of $\overline{\cg}\rlbicross \underline{\cg},$ and so $\{f^i,x^i\}$ is the dual basis for $\underline{\cg^*}\lrbicross \overline{\cg^*}.$ Denote by $\{\widetilde{\partial_i},D_i\}$ the left-invariant vector fields on $\overline{G}\rlbicross\underline{\cg}$ generated by $\{e_i, E_i\},$ and denote by $\{\widetilde{\omega^i},\widetilde{\extd x^i}\}$ the corresponding dual basis of left-invariant $1$-forms. By construction, when viewing any $f\in C^\infty(G),\ \phi\in\cg^*\subset C^\infty(\cg)$ as functions on tangent bundle, we know
\begin{gather*}
\widetilde{\partial_i} f=\partial_i f,\quad \widetilde{\partial_i}\phi=\ad^*_{e_i}\phi,\quad
D_i f=0,\quad D_i\phi= \frac{\partial}{\partial x^i}\phi.
\end{gather*} This implies 
\[\widetilde{\partial_i}=\partial_i+\sum_j(\ad^*_{e_i}x^j) \frac{\partial}{\partial x^j},\quad D_i=\frac{\partial}{\partial x^i},\quad \widetilde{\omega^i}=\omega^i,\quad \widetilde{\extd x^i}=\extd x^i-\sum_k(\ad^*_{e_k} x^i)\omega^k.\]

Let $\widetilde\circ$ be the pre-Lie structure of $\underline{\cg^*}\lcross \cg^*$ constructed by (\ref{circtilde-T}) in terms of $\ast,\circ$ in the setting of Theorem~\ref{prelie-T}. The Poisson-compatible left-covariant flat preconnection on the tangent bundle then is, for any function $a$
\begin{gather*}
\nabla_{\widehat{a}}\,{\omega^j}=\sum_{i,k}\widetilde{\partial_i} a\<f^i\ast f^j,e_k\>{\omega^k}+\sum_{i,k}D_i a\<[x^i,f^j]_{\cg^*},e_k\>{\omega^k},\\
\nabla_{\widehat{a}}\,\widetilde{\extd x^j}=\sum_{i,k} D_i a\<x^i\circ x^j,E_k\>\widetilde{\extd x^k}.
\end{gather*}
If we write
\begin{gather*}
f^i\ast f^j=\sum_k a^{ij}_k f^k,\quad x^i\circ x^j=\sum_k b^{ij}_k x^k,\\
[x^i, f^j]_{\cg^*}=\sum_k \<[x^i,f^j]_{\cg^*},e_k\> f^k=\sum_{s,k}d^{sj}_k\<x^i,e_s\>f^k,
\end{gather*}
where $[f^i,f^j]_{\cg^*}=d^{ij}_k f^k.$ Then, in terms of these structure coefficients, the left covariant preconnection on the tangent bundle $\overline{G}\rlbicross \underline{\cg}$ is
\begin{gather*}
\nabla_{\widehat{f}}\, {\omega^j}=\sum_{i,k} a^{ij}_k(\partial_i f) {\omega^k},\quad \nabla_{\widehat{f}}\, \widetilde{\extd x^j}=0;\\
\nabla_{\widehat{\phi}}\,{\omega^j}=\sum_{i,k} \left(a^{ij}_k\ad^*_{e_i}\phi+ \sum_s d^{sj}_k\<x^i,e_s\>(\frac{\partial \phi}{\partial x^i})\right)\,{\omega^k},\quad \nabla_{\widehat{\phi}}\,\widetilde{\extd x^j}=\sum_{i,k}b^{ij}_k(\frac{\partial \phi}{\partial x^i})\, \widetilde{\extd x^k},
\end{gather*}
for any $f\in C^\infty(G),\ \phi\in\cg^*\subset C^\infty(\cg).$

This result applies, for example, to tell us that we have a left-covariant differential structure on quantum groups such as $\C[\overline G]\lrbicross U_\lambda(\cg^*)$ at least to lowest order in deformation. In the special case when the product $\ast$ is zero, there is a natural differential calculus not only at lowest order. Under the notations above, we have:

\begin{proposition}\label{propC(G)U(g*)}
Let $G$ be a finite dimensional connected and simply connected Poisson-Lie group with Lie algebra $\cg.$ If the dual Lie algebra $\cg^*$ admits a pre-Lie structure $\circ:\cg^*\tens\cg^*\to\cg^*$ with respect to its Lie bracket ($[\ ,\ ]_{\cg^*}$ determined by $\delta_\cg$), then the bicrossproduct $\C[\overline{G}]\lrbicross U_\lambda(\cg^*)$ (if it exists) admits a left covariant differential calculus $\Omega^1=(\C[\overline{G}]\lrbicross U_\lambda(\cg^*))\rcross\Lambda^1$ with left invariant $1$-forms $\Lambda^1$ spanned by basis $\{\omega^i,\widetilde{\extd x^i}\}$, where the commutation relations and the derivatives are given by 
\begin{gather*}
	[f,{\omega^i}]=0,\quad [f,\widetilde{\extd x^i}]=0,\quad [x^i,{\omega^j}]=\sum_k \lambda \<[{x^i},f^j]_{\cg^*},e_k\> {\omega^k},\quad [x^i,\widetilde{\extd x^j}]=\lambda\widetilde{\extd (x^i\circ x^j)},\\
	{\extd} f=\sum_j (\partial_j f)\omega^j,\quad {\extd} x^i=\widetilde{\extd x^i}+ \sum_j(\ad^*_{e_j}x^i)\omega^j
\end{gather*}
for any $f\in \C[\overline{G}].$ This first order differential calculus extends to $C^\infty(G)\lcross U_\lambda(\cg^*)$ if one is only interested in algebra and its calculus.
\end{proposition}
\proof It is easy to see that we have a bimodule $\Omega^1$. As the notation indicates\cite{MT}, the left action on $\Omega^1$ is the product of the bicrossproduct quantum group on itself while the right action is the tensor product of the right action of the bicrossproduct on itself and a right action on $\Lambda^1$. The right action of $\C[G]$ here is trivial, namely $\omega^j\ra f=f(e)\omega^j,\ \widetilde{\extd x^j}\ra f=f(e)\widetilde{\extd x^j}$, the right actions of $x^i$ are clear from the commutation relations and given by (summation notation omitted)
\[\omega^j\ra x^i=-\lambda \<[x^i,f^j]_{\cg^*},e_k\>\omega^k=-\lambda d^{sj}_k\<x^i,e_s\>\omega^k,\quad
\widetilde{\extd x^j}\ra x^i=-\lambda\widetilde{\extd (x^i\circ x^j)}=-\lambda b^{ij}_k\widetilde{\extd x^k}.\]
One can check that these fit together to a right action of the bicrossproduct quantum group by using the Jacobi identity of $\cg^*,$ the pre-Lie identity on $\circ$, and the fact that $(\widetilde{x^i} f)(e)=\widetilde{x^i}_ef=0$ by (\ref{phitilde}).

We check that the Leibniz rule holds. The conditions $\extd [f,h]=0$ and $\extd [x^i, x^j]=\lambda\extd [x^i,x^j]_{\cg^*}$ are easy to check, so we omit these. It remains to check that 
\begin{equation}\label{cross}
\extd [x^i,f]=\lambda\extd (\widetilde{x^i} f),\quad \forall\,f\in \C[\overline{G}].
\end{equation}
The right hand side of (\ref{cross}) is \[\lambda\extd (\widetilde{x^i}f)=\lambda\partial_j(\widetilde{x^i}f)\omega^j,\]
while the left hand side of (\ref{cross}) is 
\begin{align*}
	\extd [x^i, f]&=\extd (x^i f- x^i f)=[\extd x^i,f]+[x^i,\extd f]\\
	&=[\widetilde{\extd x^i}+ (\ad^*_{e_j}x^i)\omega^j,f]+[x^i,(\partial_j f)\omega^j]\\
	&=0+[\ad^*_{e_j}x^i, f]\omega^j+[x^i,\partial_jf]\omega^j+(\partial_k f)[x^i,\omega^k]\\
	&=[\ad^*_{e_j}x^i, f]\omega^j+[x^i,\partial_jf]\omega^j+\lambda(\partial_k f)\<[x^i,f^k]_{\cg^*},e_j\>\omega^j\\
	&=\left([\ad^*_{e_j}x^i, f]+[x^i,\partial_jf]+\lambda
	\<[x^i,f^k]_{\cg^*},e_j\>(\partial_k f)\right)\omega^j\\
	&=\lambda\left(\widetilde{\ad^*_{e_j}x^i}f+\widetilde{x^i}(\partial_jf)+
	\<[x^i,f^k]_{\cg^*},e_j\>(\partial_k f)\right)\omega^j.
\end{align*}
It suffices to show that $\partial_j(\widetilde{x^i}f)=\widetilde{\ad^*_{e_j}x^i}f+\widetilde{x^i}(\partial_jf)+\<[x^i,f^k]_{\cg^*},e_j\>(\partial_k f),$
namely 
\[[\partial_j,\widetilde{x^i}]=\widetilde{\ad^*_{e_j}x^i}+
\<[x^i,f^k]_{\cg^*},e_j\>\partial_k.\]
Recall that in double cross sum $\cg^*\dcross\cg,$ for any $e_j\in\cg,\,x^i\in\cg^*$
\[[e_j, x^i]=e_j\ra x^i+e_j\la x^i=\<[x^i,f^k]_{\cg^*},e_j\>e_k+\ad^*_{e_j}x^i.\]
Therefore the condition left to check is nothing but the Lie bracket of elements $e_j,x^i$ viewing as infinitesimal action of $\cg^*\dcross\cg$ on $C^\infty(G)$ as explained in the general theory of double cross sum in Section~5.1.
\endproof

Now we compute the left covariant first order differential calculus on the bicrossproduct quantum group $\C[SU_2]\lrbicross U_\lambda(su^*_2)$ constructed in Example~\ref{TSU_2} in detail. 
\begin{example} As in Example~\ref{CqSU2}, the classical connected left-covariant calculus on $\C[SU_2]$ has basis
of left-invariant $1$-forms 
 \[\omega^0=d\extd a-b\extd c=c\extd b-a\extd d,\quad \omega^+=d\extd b-b\extd d,\quad \omega^-=a\extd c-c\extd a\]
(corresponding to the Chevalley basis $\{H,X_\pm\}$ of $su_2$) with exterior derivative
\begin{gather*}
\extd a= a \omega^0+b \omega^-,\quad \extd b=a \omega^+-b \omega^0,\\
\extd c= c \omega^0+d \omega^-,\quad \extd d=c \omega^+-d \omega^0.
\end{gather*}
Let $\circ:su^*_2\tens su^*_2\to su^*_2$ be a left pre-Lie algebra structure of $su^*_2$ with respect to the Lie bracket $[x^1,x^2]=0,$ $[x^i,x^3]=x^i,$ for $i=1,2$. Let $\{\widetilde{\extd x^1},\widetilde{\extd x^2},\widetilde{\extd x^3}\}$ complete the basis of left invariant 1-forms on the tangent bundle as explained above. According to Proposition~\ref{propC(G)U(g*)}, this defines a $6$D connected left-covariant  differential calculus over the bicrossproduct  $\C[SU_2]\lrbicross U_\lambda(su_2^*)$ with  commutation relations and exterior derivative given by
\begin{gather*}
[\mathbf{t},{\omega^l}]=0,\ \forall\,l\in\{0,\pm\},\quad [\mathbf{t},\widetilde{\extd x^i}]=0,\quad[x^i,\widetilde{\extd x^j}]=\lambda\widetilde{\extd (x^i\circ x^j)},\ \forall\,i,j\in\{1,2,3\},\\
[x^1,\omega^0]=\frac{\lambda}{2}(\omega^++\omega^-),\quad [ x^1, \omega^+]=0,\quad [x^1,\omega^-]=0,\nonumber\\
[x^2,\omega^0]=\frac{\imath \lambda}{2}(\omega^+-\omega^-),\quad [x^2,\omega^+]=0,\quad [x^2,\omega^-]=0,\\
[x^3,\omega^0]=0,\quad [x^3,\omega^+]=-\lambda\omega^+,\quad [x^3,\omega^-]=-\lambda\omega^-,\nonumber\\
\extd \begin{pmatrix} a & b\\ c & d \end{pmatrix}
=\begin{pmatrix} a & b\\ c & d \end{pmatrix} \begin{pmatrix} \omega^0 & \omega^+\\ \omega^- & -\omega^0 \end{pmatrix},\\
\extd x^1=\widetilde{\extd x^1}+2\imath x^2\omega^0+x^3\omega^+-x^3\omega^-,\quad
\extd x^2=\widetilde{\extd x^2}-2\imath x^1\omega^0+\imath x^3\omega^++\imath x^3\omega^-,\\
\extd x^3=\widetilde{\extd x^3}-(x^1+\imath x^2)\omega^++(x^1-\imath x^2)\omega^-.
\end{gather*}
\end{example}
\proof The commutation relations and derivative are computed from the formulae provided in Proposition~\ref{propC(G)U(g*)}. It is useful to also provide an independent more algebraic proof of the example from  \cite[Theorem 2.5]{MT}, where left-covariant first order differential calculi $\Omega^1$ over a Hopf algebra $H$ are constructed from pairs $(\Lambda^1,\omega)$ where $\Lambda^1$ is a right $H$-module and $\omega:H^+\to\Lambda^1$ is surjective right $H$-module map. Given such a pair, the commutation relation and derivative are given by $[h,v]=hv-h\o v\ra h\t$ and $\extd h=h\o\tens\omega\pi(h\t)$ for any $h\in H,\,v\in\Lambda^1,$ where $\Delta h=h\o\tens h\t$ denotes the coproduct of $H$ and $\pi=\id-\epsilon_H.$

Firstly, the classical calculus on $A:=\C[SU_2]$ corresponds to a pair $(\Lambda^1_A,\omega_{_{A}})$ with $\Lambda^1_A=\mathrm{span}\{\omega^0,\omega^\pm\},$ where the right $\C[SU_2]$-action on $\Lambda^1_A$ and the right $\C[SU_2]$-module surjective map $\omega_{_{A}}:\C[SU_2]^+\to \Lambda^1_A$ are given by
\[\omega^j\ra \mathbf{t}=\epsilon(\mathbf{t})\omega^j,\ j\in\{0,\pm\},\quad\omega_{_{A}}(\mathbf{t}-I_2)=\omega_{_{A}}\begin{pmatrix} a-1 & b\\ c & d-1 \end{pmatrix}=\begin{pmatrix} \omega^0 & \omega^+\\ \omega^- & -\omega^0 \end{pmatrix}.\]
Meanwhile the calculus over $H:=U_\lambda(su^*_2)$ corresponds to pair $(\Lambda^1_H,\omega_{_{H}})$ with $\Lambda^1_H=\mathrm{span}\{\widetilde{\extd x^1},\widetilde{\extd x^2},\widetilde{\extd x^3}\},$ in which the right $U_\lambda(su^*_2)$-action on $\Lambda^1_H$ and the right $U_\lambda(su^*_2)$-module surjective map $\omega_{_{H}}: U_\lambda(su^*_2)^+\to \Lambda^1_H$ are given by
\[\widetilde{\extd x^j}\ra x^i=-\lambda\widetilde{\extd (x^i\circ x^j)},\quad \omega_{_{H}}(x^i)=\widetilde{\extd x^i},\quad\forall\, i,j\in\{1,2,3\}.\]

Next we construct a pair over $\widetilde H=A\lrbicross H$ with direct sum $\Lambda^1=\Lambda^1_A\oplus\Lambda^1_H$. 
First, it is clear that $\Lambda^1_H$ is a right $\widetilde{H}$-module with trivial $A$-action
$\extd x^j\ra \mathbf{t}=\epsilon(\mathbf{t})\extd x^j,$
One can see this more generally as $v\ra ((h\o\la a)h\t)=\epsilon(a) v\ra h=(v\ra h)\ra a=v\ra (ha).$
Next, we define a right $U_\lambda(su^*_2)$-action on $\Lambda^1_A$ by the Lie bracket of $su^*_2$ viewing $\{\omega^0,\omega^\pm\}$ as $\{\phi,\psi_\pm\}$ (the dual basis to $\{H,X_\pm\}$), where $\{x^1=\psi_++\psi_-, x^2=\imath(\psi_+-\psi_-), x^3=2\phi\}$ is the basis for the half-real form $su_2^*$ of $sl_2^*,$ namely
\begin{gather}
\omega^0\ra x^1=-\frac{\lambda}{2}(\omega^++\omega^-),\quad \omega^+\ra x^1=0,\quad \omega^-\ra x^1=0,\nonumber\\
\omega^0\ra x^2=-\frac{\imath \lambda}{2}(\omega^+-\omega^-),\quad \omega^+\ra x^2=0,\quad \omega^-\ra x^2=0,\label{HonW}\\
\omega^0\ra x^3=0,\quad \omega^+\ra x^3=\lambda\omega^+,\quad \omega^-\ra x^3=\lambda\omega^-.\nonumber
\end{gather} 
This $H$-action commutes with the original trivial $A$-action on $\Lambda^1_A,$ hence $\Lambda^1_A$ also becomes a right $\widetilde{H}$-module, so does $\Lambda^1_A\oplus\Lambda^1_H.$

We then define the map $\omega:\widetilde H^+\to \Lambda^1_A\oplus \Lambda^1_H$ on generators by
\[\omega(\mathbf{t}-I_2)=\omega_{_{A}}(\mathbf{t}-I_2)=\begin{pmatrix} \omega^0 & \omega^+\\ \omega^- & -\omega^0 \end{pmatrix},\quad \omega(x^i)=\omega_{_{H}}(x^i)=\widetilde{\extd x^i},\quad\forall\, i\in\{1,2,3\}.\]
This extends to the whole $\widetilde H^+$ as a right $\widetilde{H}$-module map. To see such $\omega$ is well-defined, it suffices to check
\[\omega(x^i\mathbf{t}-\mathbf{t}x^i)=\omega([x^i,\mathbf{t}]),\quad\forall\, i\in\{1,2,3\},\]
where $[x^i,\mathbf{t}]$ are cross relations (\ref{TSU_2-action}) computed in Example~\ref{TSU_2}. On the one hand, $\omega(x^i\mathbf{t}-\mathbf{t}x^i)=\omega(x^i\mathbf{t}-(\mathbf{t}-I_2)x^i-x^i I_2)=\omega_{_{H}}(x^i)\ra\mathbf{t}-\omega_{_{A}}(\mathbf{t}-I_2)\ra x^i-\omega_{_{H}}(x^i)I_2=-\omega_{_{A}}(\mathbf{t}-I_2)\ra x^i,$ namely
\begin{equation}\label{HonW1}
\omega(x^i\mathbf{t}-\mathbf{t} x^i)=-\begin{pmatrix} \omega^0 & \omega^+\\ \omega^- & -\omega^0 \end{pmatrix}\ra x^i.
\end{equation}
Since
\begin{align*}
[x^1,\mathbf{t}]&=-\lambda bc\,\mathbf{t}e_2+\frac{\lambda}{2}\mathbf{t}\,\mathrm{diag}(ac,-bd)+\frac{\lambda}{2}\mathrm{diag}(b,-c)\\
&=-\lambda bc\,\mathbf{t}e_2+\frac{\lambda}{2}(\mathbf{t}-I_2)\,\mathrm{diag}(ac,-bd)+\frac{\lambda}{2}\mathrm{diag}(ca,-bd)+\frac{\lambda}{2}\mathrm{diag}(b,-c),
\end{align*}
we know 
\begin{align*}
\omega([x^1,\mathbf{t}])&=-\lambda \omega(b)\epsilon(c\,\mathbf{t}e_2)+\frac{\lambda}{2}\omega((\mathbf{t}-I_2))\epsilon(\mathrm{diag}(ac,-bd))\\
&+\frac{\lambda}{2}\mathrm{diag}(\omega(c)\ra a,-\omega(b)\ra d)+\frac{\lambda}{2}\mathrm{diag}(\omega(b),-\omega(c))\\
&=\frac{\lambda}{2}\mathrm{diag}(\omega^++\omega^-,-\omega^+-\omega^-),
\end{align*} as $\epsilon(\mathbf{t})=I_2.$
Likewise, we have
\begin{gather*}
\omega([x^1,\mathbf{t}])=\frac{\lambda}{2}\begin{pmatrix} \omega^++\omega^- & 0\\ 0 & -\omega^+-\omega^- \end{pmatrix},\quad \omega([x^2,\mathbf{t}])=\frac{\imath \lambda}{2}\begin{pmatrix} \omega^+-\omega^- & 0\\ 0 & -\omega^++\omega^- \end{pmatrix},\\
\omega([x^3,\mathbf{t}])=\lambda \begin{pmatrix} 0 & -\omega^+\\ -\omega^- & 0 \end{pmatrix}.
\end{gather*}
Compare with (\ref{HonW1}), we see that $\omega(x^i\mathbf{t}-\mathbf{t}x^i)=\omega([x^i,\mathbf{t}])$ holds for each $i=1,2,3$ if and only if the right $H$-action on $W$ is the one defined by (\ref{HonW}). From the coproduct of $x^i$ given in Example~\ref{TSU_2}, we know ${\extd} x^i=\widetilde{\extd x^i}+\frac{1}{2}x^k\omega(\pi(\mathrm{Tr}(\mathbf{t}\sigma_i \mathbf{t}^{-1}\sigma_k)))$. This give rise to the formulae for derivatives on $x^i$ as displayed.  
\endproof

We now analyse when a Poisson-compatible left-covariant flat preconnection is bicovariant.
\begin{lemma}\label{prelie-T-bi}
Let $\cg$ be in the setting of Theorem~\ref{prelie-T}. The pre-Lie structure $\widetilde\circ$ given by (\ref{circtilde-T}) of $\underline{\cg^*}\lrbicross \overline{\cg^*}$ obeys the corresponding (\ref{Xi-bi}) if and only if the following holds
\begin{gather}
\delta_{\cg^*}(f\circ g)=0,\quad f\o\tens [f\t,g]_{\cg^*}=0,\\
f\o\circ g\tens f\t=-f\circ g\o\tens g\t,\\
\delta_{\cg^*}(\phi\ast\psi)=0,\quad \phi\ast f\o\tens f\t=0
\end{gather} for any $\phi,\psi\in\underline{\cg^*},\,f,g\in \overline{\cg^*}.$
\end{lemma}
\proof
Since (\ref{Xi-bi}) is bilinear on entries, it suffices to show that $\widetilde\circ$ obeys (\ref{Xi-bi}) on any pair of elements $(\phi,\psi)$, $(\phi,f)$, $(f,\phi)$ and $(f,g)$ if and only if all the displayed identities holds for any $\phi,\psi\in\underline{\cg^*},\,f,g\in \overline{\cg^*}.$

Firstly, for any $f\in\overline{\cg^*},\,\phi\in\underline{\cg^*},$ 
(\ref{Xi-bi}) on $\widetilde\circ$ reduces to
\[\delta_{\cg^*}[f,\phi]_{\cg^*}-[f,\phi\o]\tens\phi\t-\phi\o\tens[f,\phi\t]=\underline{f\o}\ast\phi\tens\overline{f\t}+[\overline{f\o},\phi]\tens f\t.\]
The only term in the above identity not lying in $\underline{\cg^*}\tens\underline{\cg^*}$ is $\underline{f\o}\ast\phi\tens\overline{f\t}$ hence equals to zero. Note that $\delta_{\cg^*}$ is a $1$-cocycle, the rest of terms implies $f\o\tens [f\t,\phi]_{\cg^*}=0.$ Change the role of $f$ and $\phi$ in (\ref{Xi-bi}) implies $\phi\ast\underline{f\o}\tens\overline{f\t}=0,$ which is the same thing we just get.

Next, for any $f,g\in\overline{\cg^*},$ the condition (\ref{Xi-bi}) on $\widetilde\circ$ requires
\begin{equation*}
\begin{split}
\underline{(f\circ g)\o}\tens\overline{(f\circ g)\t}+\overline{(f\circ g)\o}\tens \underline{(f\circ g)\t}-[f,\underline{g\t}]_{\cg^*}\tens \overline{g\t}-f\circ\overline{g\o}\tens \underline{g\t}\\
-\underline{g\o}\tens f\circ\overline{g\t}-\overline{g\o}\tens [f,\underline{g\t}]_{\cg^*}\\
=[\underline{f\o},g]_{\cg^*}\tens\overline{f\t}+\overline{f\o}\circ g\tens \underline{f\t}-\underline{g\o}\tens\overline{g\t}\circ f.
\end{split}
\end{equation*}
The terms in the above identity lying in $\overline{\cg^*}\tens\underline{\cg^*}$ is exactly the condition (\ref{Xi-bi}) on pre-Lie structure $\circ$ for $\overline{\cg^*}.$ Regarding this and forgetting all the lines, the rest terms lying in $\underline{\cg^*}\tens\overline{\cg^*}$ reduces to $g\o\circ f\tens g\t+g\circ f\o\tens f\t=0,$ which is equivalent to \[f\o\circ g\tens f\t+f\circ g\o\tens g\t=0,\quad \forall\,f,g\in \cg^*.\] Combining above wih $f\o\tens [f\t,\phi]_{\cg^*}=0,$ we know the condition (\ref{Xi-bi}) on $\circ$ reduce to $\delta_{\cg^*}(f\circ g)=0.$

Finally, for any $\phi,\psi\in\underline{\cg^*},$ the condition (\ref{Xi-bi}) on $\widetilde\circ$ reduces to (\ref{Xi-bi}) on $\ast$ for $\underline{\cg^*}.$ Since $\ast$ is commutative, this eventually becomes
\[(\phi\ast\psi)\o\tens(\phi\ast\psi)\t=\phi\ast\psi\o\tens\psi\t+\phi\o\ast\psi\tens\phi\t.\] Since $\phi\ast f\o\tens f\t=0,$ this reduces to $\delta_{\cg^*}(\phi\ast \psi)=0.$ This finishes our proof.
\endproof

The conditions in Lemma~\ref{prelie-T-bi} all hold when the Lie bracket of $\cg$ (or the Lie cobracket of $\cg^*$) vanishes. Putting these results together we have:

\begin{proposition}\label{tangentthm}
Let $G$ be a finite dimensional connected and simply connected Poisson-Lie group with Lie bialgebra $\cg.$ Assume that $(\cg^*,[\ ,\ ]_{\cg^*})$ obeys the conditions in Theorem~\ref{prelie-T} and Lemma~\ref{prelie-T-bi}. Then the tangent bundle $\overline{G}\rlbicross\underline{\cg}$ in Lemma~\ref{tangent} admits a Poisson-compatible bicovariant flat preconnection.
\end{proposition}

\begin{example}  In the setting of Example~\ref{Tm^*}, we already know from Corollary~\ref{preliecorol} that the abelian Poisson-Lie group $\R^n\lcocross\cm^*$ admits a Poisson-compatible left-covariant (bicovariant) flat preconnection if and only if  $(\overline{\cm^*}\lcocross\cm^*)^*=\underline{\cm}\lcross_\ad\cm$ admits a pre-Lie structure. 

From Corollary~\ref{semipre-m}, we know such pre-Lie structure $\widetilde\circ$ exists and is given by $(x,\xi)\widetilde\circ(y,\eta)=(x\cdot y+[\xi,y]_\cm,\xi\circ\eta)$ if we assume $(\cm,\cdot,[\ ,\ ]_\cm)$ to be a finite dimensional (not necessarily unital) Poisson algebra such that $(\cm, [\ ,\ ]_\cm)$ admits a pre-Lie structure $\circ:\cm\tens\cm\to\cm.$ Then the corresponding preconnection is
\[\nabla_{\widehat{(x,\xi)}}\extd (y,\eta)=\extd (x\cdot y+[\xi,y]_\cm,\xi\circ\eta)\]
for any $x,y\in\underline{\cm},\,\xi,\eta\in\cm.$ 

In fact this extends to all orders.
Under the assumptions above, according to Proposition~\ref{envel}, the noncommutative algebra $U_\lambda(\underline{\cm}\lcross_\ad\cm)=S(\underline{\cm})\lcross U_\lambda(\cm),$ or the cross product of algebras $\C[\R^n]\lcross U(\cm)$ (as quatisation of $C^\infty(\R^n\lcocross\cm^*)$), admits a connected bicovariant differential graded algebra \[\Omega(U_\lambda(\underline{\cm}\lcross_\ad\cm))=(S(\underline{\cm})\lcross U_\lambda(\cm))\rcross \Lambda(\underline{\cm}\lcross_\ad\cm)\]
as quantised differential graded algebra. Note that $\extd (x,\xi)=1\tens (x,\xi)\in 1\tens\Lambda^1.$ The commutation relations on generators are
\begin{gather*}
[\xi,\eta]=\lambda[\xi,\eta]_\cm,\quad [x,y]=0,\quad [\xi,x]=\lambda[\xi,x]_\cm,\\
[x,\extd y]=\lambda\extd(x\cdot y),\quad [\xi,\extd x]=\lambda\extd [\xi,x]_{\cm},\quad [\xi,\extd\eta]=\lambda \extd(\xi\circ\eta)
\end{gather*}
for any $x,y\in\underline{\cm},\,\xi,\eta\in\cm.$
\end{example}

\section{Quantisation of cotangent bundle $T^*G=\underline{\cg^*}\lcross G$}
In this section, we focus on quantisation of cotangent bundle $T^*G$ of a Poisson-Lie group $G.$ We aim to construct preconnections on $T^*G.$

As a Lie group, the cotangent bundle $T^*G$ can be identified with the semidirect product of Lie groups $\underline{\cg^*}\lcross G$ with product given by
\begin{equation*}
(\phi,g)(\psi,h)=(\phi+\Ad^*(g)(\psi),gh)
\end{equation*}
for any $g,h\in G,\,\phi,\psi\in\cg^*.$ As before, $\underline{\cg^*}$ is $\cg^*$ but viewed as abelian Lie group under addition.
In particular, $(\phi,g)^{-1}=(-\Ad^*(g^{-1})(\phi),g^{-1})$ and $(0,g)(\phi,e)(0,g)^{-1}=(\Ad^*(g)\phi,e).$
Here $\Ad^*$ is the coadjoint action of $G$ on the dual of its Lie algebra. The Lie algebra of $T^*G$ is then identified with the semidirect sum of Lie algebras $\underline{\cg^*}\lcross \cg,$ where the Lie bracket of $\underline{\cg^*}\lcross \cg$ is given by
\begin{equation}\label{semi-lie}
[(\phi,x),(\psi,y)]=(\ad^*_x\psi-\ad^*_y\phi,[x,y]_{\cg})
\end{equation}
for any $\phi,\psi\in\underline{\cg^*},\,x,y\in\cg.$ 
Here $\underline{\cg^*}$ is $\cg^*$ viewed as abelian Lie algebra and $\ad^*$ denotes the usual left coadjoint action of $\cg$ on $\cg^*$ (or $\underline{\cg}^*$). 

The strategy to build Poisson-Lie structure on cotangent bundle here is to construct Lie bialgebra structures on $\underline{\cg^*}\lcross\cg$ via bosonization of Lie bialgebras. Then we can exponentiate the obtained Lie cobracket of $\underline{\cg^*}\lcross\cg$ to a Poisson-Lie structure on $\underline{\cg^*}\lcross G.$ We can always do this, as we work in the nice case when Lie group is connected and simply connected.

\subsection{Lie bialgebra structures on $\underline{\cg^*}\lcross\cg$ via bosonization}
Let ${}^{\cg}_{\cg}\Mcal$ denote the monoidal category of left Lie $\cg$-crossed modules. A \textit{braided-Lie bialgebra} $\mathfrak{b}\in {}^{\cg}_{\cg}\Mcal$ is  $(\mathfrak{b},[\ ,\ ]_{\mathfrak{b}},\delta_{\mathfrak{b}},\la,\beta)$ given by a $\cg$-crossed module $(\mathfrak{b},\la,\beta)$ that is both Lie algebra $(\mathfrak{b},[\ ,\ ]_{\mathfrak{b}})$ and Lie coalgebra $(\mathfrak{b},\delta_{\mathfrak{b}})$ live in ${}^{\cg}_{\cg}\Mcal,$ and the infinitesimal braiding $\Psi:\mathfrak{b}\tens \mathfrak{b}\to \mathfrak{b}\tens \mathfrak{b}$ obeying 
$\Psi(x,y)=\ad_x\delta_{\mathfrak{b}} y-\ad_y\delta_{\mathfrak{b}} x-\delta_{\mathfrak{b}}([x,y]_{\mathfrak{b}})$ for any $x,y\in {\mathfrak{b}}.$  
If $\mathfrak{b}$ is a braided-Lie bialgebra in ${}^{\cg}_{\cg}\Mcal,$ then the bisum $\mathfrak{b}\lbiprod \cg$ with semidirect Lie bracket/cobracket is a Lie bialgebra~\cite{Ma:blie}.  

For our purpose, a straightforward solution is to ask for $\underline{\cg^*}=(\cg^*, [\ ,\ ]=0,\delta_{\cg^*},\ad^*,\alpha)$ being a braided-Lie algebra in ${}^{\cg}_{\cg}\Mcal$ for some left $\cg$-coaction $\alpha$ on $\cg^*.$ 

\begin{lemma}\label{blie}
Let $\cg$ be a finite-dimensional Lie bialgebra and suppose there is a linear map $\Xi:\cg^*\tens\cg^* \to \cg^*$ such that (\ref{compatible}) holds. Then $\underline{\cg^*}=(\cg^*,[\ ,\ ]=0,\delta_{\cg^*},\ad^*,\alpha)$ (a variation of $\cg^*$) is a braided-Lie bialgebra in ${}^{\cg}_{\cg}\Mcal$ if and only if   $\Xi$ is a pre-Lie structure on $\cg^*$ satisfying that $\Xi$ is covariant under Lie cobracket $\delta_{\cg^*}$ in the sense of
\begin{equation}\label{Xi-ass}
\Xi(\phi,\psi)\o\tens \Xi(\phi,\psi)\t=\Xi(\phi, \psi\o)\tens\psi\t+\psi\o\tens\Xi(\phi,\psi\t),
\end{equation} 
and 
\begin{equation}\label{Xi-con}
\Xi(\phi\o,\psi)\tens\phi\t=\psi\o\tens\Xi(\psi\t,\phi)
\end{equation}
for any $\phi,\psi\in\cg^*$. Here the left $\cg$-coaction $\alpha$ and the left pre-Lie product $\Xi$ of $\cg^*$ are mutually determined via
\begin{equation}\label{alpha-Xi}
\<\alpha(\phi),\psi\tens x\>=-\Xi(\psi,\phi)(x),
\end{equation}
for any $\phi,\psi\in\cg^*,\,x\in\cg.$ In this case, the bisum $\underline{\cg^*}\lbiprod \cg$ is a Lie bialgebra with Lie bracket given by (\ref{semi-lie}) and Lie cobracket given by
\begin{equation}\label{semi-lie-co}
\delta(\phi,X)=\delta_{\cg}X+\delta_{\cg^*}\phi+(\id-\tau)\alpha(\phi)
\end{equation}
for any $\phi\in\cg^*,X\in\cg.$
\end{lemma}
\proof 
Since Lie bracket is zero, by definition, the question amounts to find left $\cg$-coaction $\alpha$ on $\cg^*$ such that 1) $(\ad^*,\alpha)$ makes $\underline{\cg^*}$ into a left $\cg$-crossed module; 2) $\delta_{\cg^*}$ is a left $\cg$-comodule map under $\alpha$; 3) The infinitesimal braiding $\Psi$ on $\underline{\cg^*}$ is trivial, i.e.,
\begin{equation}\label{Psi}
\Psi(\phi,\psi)=\ad^*_{\psi\bo}\phi\tens\psi\bt-\ad^*_{\phi\bo}\psi\tens\phi\bt
-\psi\bt\tens\ad^*_{\psi\bo}\phi+\phi\bt\tens\ad^*_{\phi\bo}\psi
\end{equation}
is zero for any $\phi,\psi\in\underline{\cg^*},$ where we denote $\alpha(\phi)=\phi\bo\tens\phi\bt.$

Clearly, $\alpha$ is a left $\cg$-coaction on $\cg^*$ if and only if $\Xi$ defines a left $\cg^*$ action on itself, since $\alpha$ and $\Xi$ are adjoint to each other by (\ref{alpha-Xi}), thus if and only if $\Xi$ is left pre-Lie structure, due to (\ref{compatible}). Next, the condition that the Lie cobracket $\delta_{\cg^*}$ is a left $\cg$-comodule map under $\alpha$ means $\delta_{\cg^*}$ is a right $\cg^*$-module map under $-\Xi.$ This is exactly the assumption (\ref{Xi-ass}) on $\Xi.$ In this case, the cross condition (\ref{g-crossv}) or (\ref{Xi-bi}) (using compatibility) for making $\cg^*$ a left $\cg$-crossed module under $(\ad^*,\alpha)$ becomes (\ref{Xi-con}). 

It suffices to show that the infinitesimal braiding $\Psi$ on $\cg^*$ is trivial on $\underline{\cg}^*.$ By construction, $\<\alpha(\phi),\varphi\tens x\>=-\Xi(\varphi,\phi)(x),$ so $\ad^*_{\psi\bo}\phi\tens\psi\bt=\phi\t\tens\Xi(\phi\o,\psi)$ where $\alpha(\phi)=\phi\bo\tens\phi\bt$ and $\delta_{\cg^*}\phi=\phi\o\tens\phi\t.$ Thus
\begin{align*}
\Psi(\phi,\psi)&=\ad^*_{\psi\bo}\phi\tens\psi\bt-\ad^*_{\phi\bo}\psi\tens\phi\bt
-\psi\bt\tens\ad^*_{\psi\bo}\phi+\phi\bt\tens\ad^*_{\phi\bo}\psi\\
&=
\Xi(\psi\o,\phi)\tens\psi\t-\psi\t\tens\Xi(\psi\o,\phi)
-\Xi(\phi\o,\psi)\tens\phi\t+\phi\t\tens\Xi(\phi\o,\psi)\\
&=\Xi(\psi\o,\phi)\tens\psi\t+\psi\o\tens\Xi(\psi\t,\phi)
-\Xi(\phi\o,\psi)\tens\phi\t-\phi\o\tens\Xi(\phi\t,\psi)\\
&=0,
\end{align*} using (\ref{Xi-con}).
This finishes our proof. 
\endproof

\begin{example}\label{kk-blie} Let $\cm$ be a pre-Lie algebra with product $\circ:\cm\tens\cm\to\cm$ and $\cg=\cm^*$ with zero Lie bracket as in Example~\ref{kk}. This meets the conditions in Lemma~\ref{blie} and we have a Lie bialgebra
$\underline{\cg^*}\lcocross\cg=\underline{\cm}\lcocross\cm^*$ with zero Lie bracket and with Lie cobracket
\[ \delta \phi=\delta_{\cm^*}\phi,\quad \delta x= (\id-\tau)\alpha(x),\quad\forall \phi\in \cm^*,\ x\in\cm\]
where $\alpha$ is given by the pre-Lie algebra structure $\circ$ on $\cm$, i.e., $\<x\tens\phi,\alpha(y)\>=-\<\phi,x\circ y\>.$ The Lie bialgebra here is the dual of semidirect sum Lie algebra $\tilde \cm=\cm^*\lcross\cm$ (viewed as Lie bialgebra with zero Lie cobracket) where $\cm$ acts on $\cm^*$ by the adjoint action of $\cm$ on $\cm$ given by $\circ$, i.e., $\<x\la\phi,y\>=-\phi(x\circ y)$, 
\[ [x,y]=[x,y]_\cm,\quad  [x,\phi]=x\la\phi,\quad [\phi,\psi]=0,\quad \forall x,y\in\cm,\ \phi,\psi\in\cm^*.\]
The Poisson bracket on $\tilde{\cm}^*=\underline{\cm}\lcocross\cm^*$ is then the Kirillov-Kostant one for $\tilde\cm$, i.e., given by this Lie bracket. \end{example}

\begin{example}\label{qt-blie}
Let $\cg$ be a quasitriangular Lie bialgebra with $r$-matrix $r=r\bo\tens r\bt\in \cg\tens\cg$ such that $r_+\la X=0$ on $\cg$. As in Example~\ref{qt-flat}, $\cg^*$ is a pre-Lie algebra with product $\Xi(\phi,\psi)=-\<\phi,r\bt\>\ad^*_{r\bo}\psi.$ Direct computation shows $\Xi$ satisfies (\ref{Xi-ass}), (\ref{Xi-con}) without any further requirement.  So $\underline{\cg^*}=(\cg^*,[\ ,\ ]=0,\delta_{\cg^*},\ad^*,\alpha)$ is a braided-Lie bialgebra in ${}^{\cg}_{\cg}\Mcal$ with $\alpha(\phi)=r\bt\tens\ad^*_{r\bo}\phi.$ Hence from Lemma~\ref{blie}, $\underline{\cg^*}\lbiprod\cg$ is a Lie bialgebra with Lie bracket given by (\ref{semi-lie}) and Lie cobracket  given by (\ref{semi-lie-co}), i.e.,
\begin{equation}
\delta(\phi,X)=\delta_{\cg}X+\delta_{\cg^*}\phi+(\id-\tau)(r\bt\tens\ad^*_{r\bo}\phi).
\end{equation}

Note that if  $\cg$ is a quasitriangular Lie bialgebra, it is shown in~\cite[Corollary 3.2, Lemma 3.4]{Ma:blie} that $(\cg^*,\delta_{\cg^*})$ is a braided-Lie algebra with Lie bracket given by 
\begin{equation*}
[\phi,\psi]=2\<\phi,r_+\bo\>\ad^*_{r_+\bt}\psi=0
\end{equation*}
in our case, so in this example $\und{\cg}^*$ in Lemma~\ref{blie} agrees with a canonical construction. \end{example}

\subsection{Poisson-Lie structures on $\underline{\cg^*}\lcross G$ induced from $\underline{\cg^*}\lbiprod\,\cg$}

Next we exponentiate our Lie bialgebra structure  $\underline{\cg^*}\lbiprod \cg$ constructed by Lemma~\ref{blie} to a Poisson-Lie structure on cotangent bundle. As usual this is done by
exponentiating $\delta$ to a group 1-cocycle $D$. 

\begin{proposition}\label{bpoisson}
Let $G$ be connected and simply connected Poisson-Lie group. If its Lie algebra $\cg$ with a given coaction $\alpha$ is in the setting of Lemma~\ref{blie}. Then $\underline{\cg^*}\lcross G$ is a Poisson-Lie group with
$$D(\phi,g)=\Ad_{\phi} D(g)+\delta_{\cg^*}\phi+(\id-\tau)(\phi\bo\tens\phi\bt-\frac{1}{2}\ad^*_{\phi\bo}\phi\tens\phi\bt),$$ where $\alpha(\phi)=\phi\bo\tens\phi\bt.$ 
\end{proposition}
\proof Because of the cocycle condition it suffices to find $D(\phi):= D(\phi,e)$ and $D(g):=D(e,g)$, then
\[ D(\phi,g)=D(\phi)+\Ad_\phi D(g),\quad \forall (\phi,g)\in \cg^*\lcross G\]
where
\[ \Ad_\phi(X)=X- \ad^*_X\phi,\quad\forall X\in \cg\subset\cg^*\lcross\cg,\quad \phi\in \cg^*.\]
We require
\[ {\extd \over\extd t}D(t\phi)=\Ad_{t\phi}(\delta\phi)\]
which we solve writing
\[ D(\phi)=\delta_{\cg^*}\phi+ Z(\phi)\]
so that
\[  {\extd \over\extd t}Z(t\phi)=\Ad_{t\phi}((\id-\tau)\circ\alpha(\phi))=(\id-\tau)\circ\alpha(\phi)-t(\id-\tau)(\ad^*_{\phi\bo}\phi\tens\phi\bt),\ \ Z(0)=0.\]
Integrating this to
\[Z(t\phi)=t(\id-\tau)\circ\alpha(\phi)-{1\over 2} t^2(\id-\tau)(\ad^*_{\phi\bo}\phi\tens\phi\bt),\]
we obtain 
\[ D(\phi)=\delta_{\cg^*}\phi+(\id-\tau)(\phi\bo\tens\phi\bt-\frac{1}{2}\ad^*_{\phi\bo}\phi\tens\phi\bt),\] where $\alpha(\phi)=\phi\bo\tens\phi\bt.$
The general case $\frac{\extd}{\extd t}|_{t=0}D(\phi+t\psi)=\Ad_{\phi}(\delta\psi)$ amounts to vanishing of the expression (\ref{Psi}) which we saw 
holds under our assumptions in the proof of Lemma~\ref{blie}. 
\endproof

\begin{example} In the setting of Example~\ref{qt-blie} with $(\cg,r)$ quasitriangular such that $r_+\la X=0$ on $X\in\cg$ we know that $\cg^*\lcross G$ is a Poisson-Lie group with
\[ D(\phi,g)=\delta_{\cg^*}\phi+\Ad_{(\phi,g)}(r)-r+2r_+\la\phi - r_+\la(\phi\tens\phi),\]
where $\la$ denotes coadjoint action $\ad^*$.
As $\alpha(\phi)=r_{21}\la\phi,$ direct computation shows that
$D(\phi)=\delta_{\cg^*}(\phi)+(\id-\tau)r_{21}\la\phi+ r_-\la(\phi\tens\phi)$. Since the differential equation for $D(g)$ is the same as one on $G$, so $D(g)=\Ad_g(r)-r$ as $\cg$ quasitriangular, we obtain the stated result. Note that
$\Ad_\phi(r)=(r\bo-r\bo\la\phi)\tens (r\bt-r\bt\la\phi)=r+r\la(\phi\tens\phi)-r\bo\la\phi\tens r\bt-r\bo\tens r\bt\la\phi.$  The differential equation $\frac{\extd}{\extd t}|_{t=0}D(\phi+t\psi)=\Ad_{\phi}(\delta\psi)$ amounts to $r_+\la(\id-\tau)(\phi\tens\psi)=0$,  which is guaranteed by $r_+\la X=0$ on $\cg.$

Note that we can view $r\in (\underline\cg^*\lbiprod\cg)^{\tens 2}$ and it will obey the CYBE and in our case $\ad_\phi(r_+)=0$ as $r_+\la\phi=0$ on $\cg^*$ under our assumptions. In this case $\underline\cg^*\lbiprod\cg$ is quasitriangular with the same $r$, with Lie cobracket
\[ \delta_r(\phi)= \ad_\phi(r)=-r\bo\la\phi\tens r\bt-r\bo\tens r\bt\la\phi=(\id-\tau)r_{21}\la\phi\]
at the Lie algebra level (differentiating the above $\Ad_{t\phi}$) and $\delta X$ as before. In our case the cobracket has
an additional $\delta_{\cg^*}\phi$ term reflected also in $D$. 
\end{example}

\subsection{Preconnections on cotangent bundle $\underline{\cg^*}\lbiprod\, G$} 

Let $\cg$ be a finite-dimensional of Lie bialgebra and its dual $\cg^*$ admits a pre-Lie structure $\Xi:\cg^*\tens\cg^*\to \cg^*$ such that (\ref{Xi-ass}) and (\ref{Xi-con}) as in the setting of Lemma~\ref{blie}. Then the dual of Lie bialgebra $\underline{\cg^*}\lbiprod \cg$ is $\overline{\cg}\lbiprod \cg^*,$whose Lie bracket is the semidirect sum $\cg\lcross\cg^*$ and Lie cobracket is the semidirect cobracket $\overline{\cg}\lcocross\cg^*,$ namely
\begin{gather*}
[x,y]=[x,y]_\cg,\quad [\phi,x]=\phi\la x,\quad [\phi,\psi]=[\phi,\psi]_{\cg^*};\\
\delta x=(\id-\tau)\beta(x),\quad \delta\phi=\delta_{\cg^*}\phi,
\end{gather*} 
for any $x,y\in\cg,\,\phi,\psi\in\cg^*.$ Here the left action and coaction of $\cg^*$ on $\cg$ are given by
\begin{equation}\label{Xi-action}
\<\phi\la x, \psi\>=-\Xi(\phi,\psi)(x),\quad\text{ and    }\quad
\<\beta(x),y\tens\phi\>=\<\phi,[x,y]\>,
\end{equation}
respectively.

Here again, we use Lemma~\ref{semipre} to construct pre-Lie algebra structures on semidirect sum $\overline{\cg}\lcross\cg^*$.

\begin{theorem}\label{prelie-C} Let $G$ be a connected and simply connected Poisson-Lie group with Lie bialgebra $\cg$.  Let $\cg^*$ admit  two pre-Lie structures $\Xi$ and $\circ$, with $\Xi$ obeying (\ref{Xi-ass}) and (\ref{Xi-con}) as in the setting of Lemma~\ref{blie}. Let $\cg$ also admit a pre-Lie structure $\ast$ such that
\begin{equation}\label{Xi-ast}
\phi\la (x\ast y)=(\phi\la x)\ast y+ x\ast (\phi\la y),
\end{equation}
where $\la$ is defined by (\ref{Xi-action}). Then the Lie algebra $\overline{\cg}\lcross\cg^*$ admits a pre-Lie structure $\widetilde{\circ}:$
\begin{equation}\label{circtilde-C}
(x,\phi)\widetilde{\circ}(y,\psi)=(x\ast y + \phi\la y, \phi\circ \psi)
\end{equation}
and then the cotangent bundle $\underline{\cg^*}\lcross G$ admits a Poisson-compatible left-covariant flat preconnection.
\end{theorem}
\proof Since $(\cg,\Xi)$ is in the setting of Lemma~\ref{blie}, the left $\cg^*$-action in semidirect sum $\overline{\cg}\lcross\cg^*$ is the one defined in (\ref{Xi-action}). The rest is immediate from Lemma~\ref{semipre} and Corollary~\ref{preliecorol}. \endproof

To construct a bicovariant preconnection, the pre-Lie structure constructed in Theorem~\ref{prelie-C} must satisfy corresponding (\ref{Xi-bi}).

\begin{proposition}\label{prelie-C-bi}
In the setting of Theorem~\ref{prelie-C}, the pre-Lie structure $\widetilde\circ$ of $\overline{\cg}\lcross\cg^*$ defined by (\ref{circtilde-C}) obeys the corresponding (\ref{Xi-bi}) if and only if $\circ$ obeys (\ref{Xi-bi}), $\ast$ is associative and the following identities hold
\begin{gather}
[x,y]\ast z=[y,z]\ast x,\label{ast}\\
((\ad^*_x\psi)\circ\phi)(y)+\Xi(\ad^*_y\phi,\psi)(x)=0,\label{circad}\\
\Xi(\phi,\psi)([x,y]_\cg)=\Xi(\phi,\ad^*_y\psi)(x)-(\phi\circ\ad^*_x\psi)(y),\label{Xiad}
\end{gather} for any $x,y,z\in\cg$ and $\phi,\psi\in\cg^*.$ The associated preconnection is then bicovariant.
\end{proposition}
\proof
Since (\ref{Xi-bi}) is bilinear, it suffices to show that (\ref{Xi-bi}) holds on any pair of elements $(x,y)$, $(x,\phi)$, $(\phi,x)$ and $(\phi,\psi)$ if and only if all the conditions and  displayed identities hold. Here we denote $\beta(x)=x^1\tens x_2\in\cg^*\tens\cg,$ so we know
\[\<x^1,y\>x_2=[x,y]_{\cg},\quad x^1\<x_2,\phi\>=-\ad^*_x\phi.\]

Firstly, for any $\phi,\psi\in\cg^*,$ the condition (\ref{Xi-bi}) holds for $\widetilde\circ$ reduces to (\ref{Xi-bi}) on pre-Lie structure $\circ$ for $\cg^*.$ 

Secondly, for any $x,y\in\overline{\cg},$ the condition (\ref{Xi-bi}) requires
\begin{equation*}
\begin{split}
(x\ast y)^1\tens (x\ast y)_2- (x\ast y)_2\tens (x\ast y)^1-x^1\la y\tens x_2+x_2\ast y\tens x^1+x\ast y_2\tens y^1\\
=y^1\tens [x,y_2]_{\cg}+y_2\tens y^1\la x.
\end{split}
\end{equation*}
The terms lying in $\overline{\cg}\tens \overline{\cg}$ on both sides should be equal, i.e., $-x^1\la y\tens x_2=y_2\tens y^1\la x,$ which is equivalent to $-\Xi(\ad^*_x\psi,\phi)(y)=\Xi(\ad^*_y\phi,\psi)(x).$ This is true from our assumption (\ref{Xi-con}) on $\Xi.$ The terms in $\overline{\cg}\tens\cg^*$ is equivalent to $[x\ast y,z]=[x,z]\ast y+x\ast[y,z],$ i.e., $\ast$ is associative. The terms in $\cg^*\tens\overline{\cg}$ is $(x\ast y)^1\tens (x\ast y)_2=y^1\tens [x,y_2]_{\cg},$ apply the first entry to $z\in\cg,$ we get $[x\ast y,z]=[x,[y,z]],$ which is equivalent to $[x,z]\ast y=[z,y]\ast x.$

Now, for any $x\in\cg,\phi\in\cg^*,$ the condition (\ref{Xi-bi}) reduces to $0=x^1\circ\phi\tens x_2-\phi\o\tens\phi\t\la x.$ Apply $y\tens\psi,$ this becomes $-\Xi(\ad^*_y\phi,\psi)(x)=((\ad^*_x\psi)\circ\phi)(y).$

Finally, for any $\phi\in\cg^*,\,x\in\overline{\cg},$ the condition (\ref{Xi-bi}) requires
\begin{equation*}
\begin{split}
(\phi\la x)^1\tens (\phi\la x)_2- (\phi\la x)_2\tens (\phi\la x)^1-\phi\circ x^1\tens x_2+\phi\la x_2\tens x^1\\
-x^1\tens\phi\la x_2+x_2\tens\phi\circ x^1
=\phi\o\la x\tens \phi\t+x_2\tens x^1\circ\phi.
\end{split}
\end{equation*}
The terms lying in $\cg^*\tens\overline{\cg}$ give $(\phi\la x)^1\tens (\phi\la x)_2-\phi\circ x^1\tens x_2-x^1\tens \phi\la x_2=0.$ Apply $y\tens\psi,$ this is equivalent to 
$-\Xi(\phi,\ad^*_y\psi)(x)+(\phi\circ \ad^*_x\psi)(y)+\Xi(\phi,\psi)([x,y]_\cg)=0.$ Apply $\psi\tens y$ to terms lying $\overline{\cg}\tens\cg^*,$ after cancelling the identity we just get, we have $((\ad^*_x\psi)\circ\phi)(y)+\Xi(\ad^*_y\phi,\psi)(x)=0.$ This finishes our proof.
\endproof

For simplicity, one can certainly choose $\Xi=\circ$ in Theorem~\ref{prelie-C} and Proposition~\ref{prelie-C-bi}.

\begin{corollary}\label{prelie-C-cor}
Let $\cg$ be finite dimensional Lie bialgebra. Assume that $\cg^*$ admits a pre-Lie structure $\Xi$ such that (\ref{Xi-ass}) and (\ref{Xi-con}). Also assume that $\cg$ admits a pre-Lie structure $\ast$ such that (\ref{Xi-ast}) where the action is defined by (\ref{Xi-action}) from $\Xi$. Then
\begin{equation*}
(x,\phi)\widetilde{\circ} (y,\psi)=(x\ast y+\phi\la y,\Xi(\phi,\psi))
\end{equation*} defines a pre-Lie structure of Lie algebra $\overline{\cg}\lcross\cg^*,$ and thus provides a Poisson-compatible left-covariant flat preconnection on cotangent bundle $\cg^*\lbiprod\, G.$ Moreover, if $\ast$ is associative and obeys (\ref{ast}), then the pre-Lie structure $\widetilde{\circ}$ obeys (\ref{Xi-bi}), thus the corresponding preconnection is bicovariant.
\end{corollary}
\proof Clearly, there is no further condition on $\circ$ in the case $\circ=\Xi$ in Theorem~\ref{prelie-C}. In bicovariant case, the further conditions on $\circ$ in Proposition~\ref{prelie-C-bi} are (\ref{Xi-bi}), (\ref{circad}) and (\ref{Xiad}). These all can be showed by the assumptions (\ref{Xi-ass}) and (\ref{Xi-con}) we already made on $\Xi$. In particular, (\ref{Xi-con}) can show (\ref{circad}) is true, and (\ref{Xi-ass}) is simply a variation of (\ref{Xiad}) when $\circ=\Xi$. The only conditions left in Proposition~\ref{prelie-C-bi} are $\ast$ is associative and (\ref{ast}). This completes our proof.
\endproof

\begin{example}\label{kk-pc}
In the case of Example~\ref{kk-blie} we know the answer: a Poisson-compatible bicovariant flat preconnection on $\tilde\cm^*=\underline{\cm}\lcocross\cm^*$ corresponds to a pre-Lie algebra structure on $\tilde\cm=\cm^*\lcross\cm$. 

Assume $\tilde{\circ}$ is such a pre-Lie structure, and also assume $\tilde{\circ}$ is such that
$\tilde\circ(\cm\tens\cm)\subseteq\cm,\,\tilde\circ(\cm^*\tens\cm^*)\subseteq\cm^*,\,\tilde\circ(\cm\tens\cm^*)\subseteq\cm^*$ and the restriction of $\tilde\circ$ on other subspace is zero. Directly from the definition of pre-Lie structure, one can show $\circ:=\tilde\circ|_{\cm\tens\cm}$ also provides a pre-Lie structures for $(\cm,[\ ,\ ]_{\cm})$, while $\ast:=\tilde\circ|_{\cm^*\tens\cm^*}$ provides a pre-Lie structure for $(\cm^*,[\ ,\ ]_{\cm^*}=0),$ thus $\ast$ is associative and (\ref{ast}) holds automatically. Besides $\la:=\tilde\circ|_{\cm\tens\cm^*}$ can be shown to be a left $\cm$-action on $\cm^*$, which is exactly the adjoint action of left $\cm$-action on $\cm$ given by pre-Lie structure $\circ$ on $\cm$. Apply $\tilde{\circ}$ to any $x\in\cm,\,\phi,\psi\in\cm^*,$ one has $x\la(\phi\ast\psi)=(x\la\phi)\ast\psi+\phi\ast(x\la\psi)$, i.e., (\ref{Xi-ast}). The analysis above shows that $\circ,\ast,\la$ corresponds to the data in Corollary~\ref{prelie-C-cor}. So this example agrees with our construction of Poisson-compatible bicovariant flat preconnection on $\underline{\cg}^*\lbiprod\cg=\underline{\cm}\lcocross\cm^*$ in case of $\cg=(\cm^*,[\ ,\ ]_{\cm^*}=0)$ in Corollary~\ref{prelie-C-cor}.

We already  know how to quantise the algebra $C^\infty(\tilde{\cm}^*)$ or $S(\tilde{\cm})$ and its differential graded algebra as in Example~\ref{kk}.
More precisely, the quantisation of $S(\tilde{\cm})$ is the noncommutative algebra $U_\lambda(\tilde{\cm})$ with relation $xy-yx=\lambda [x,y]$ for any $x,y\in\tilde{\cm},$ namely \[U_\lambda(\tilde{\cm})=U_\lambda(\cm^*\lcross\cm)=S(\cm^*)\lcross U_\lambda(\cm)\] with relation $x\phi-\phi x=\lambda x\la\phi$ for any $x\in\cm,\,\phi\in\cm^*.$
Besides, as in Example~\ref{kk} and Proposition~\ref{envel}. the preconnection on $\tilde{\cm}^*=\underline\cm\lcocross\cm^*$ is given by
\[\nabla_{\widetilde{(\phi,x)}}\extd((\psi,y))=\extd((\phi,x)\tilde{\circ} (\psi,y))=\extd(\phi\ast\psi+x\la\phi,x\circ y).\]
Thus, the quantised differential calculus is \[\Omega(U_\lambda(\tilde{\cm}))=U_\lambda(\tilde{\cm})\rcross \Lambda(\tilde{\cm})=(S(\cm^*)\lcross U_\lambda({\cm}))\rcross \Lambda(\cm^*\oplus\cm)\]
with relation 
\[[(\phi,x),\extd(\psi,y)]=\lambda\,\extd (\phi\ast\psi+x\la\phi,x\circ y)\]
for any $(\phi,x), (\psi,y)\in \tilde{\cm}\subset U_\lambda(\tilde{\cm}),$ where $\Lambda(\cm^*\oplus\cm)$ denotes the usual exterior algebra on vector space $\cm^*\oplus\cm$ and $\extd (\psi,y)=1\tens (\psi+y)\in 1\tens\Lambda.$
\end{example}

\begin{example}\label{qt-pc}
In the case of $\cg$ quasitriangular with $r_+\la x=0$ on $\cg$ as in Example~\ref{qt-blie}. According to Corollary~\ref{prelie-C-cor}, if $\cg$ admits a pre-Lie product $\ast$ such that 
\[[r\bo,x\ast y]\tens r\bt=[r\bo,x]\ast y\tens r\bt+x\ast [r\bo,y]\tens r\bt,\] i.e., corresponding (\ref{Xi-ast}), then $\overline{\cg}\lbiprod\cg^*$ in Example~\ref{qt-blie} admits a pre-Lie structure $\widetilde\circ$
\[x\widetilde\circ y=x\ast y,\quad \phi\widetilde\circ x=\phi\la x=-\<\phi,r\bt\>[r\bo, x],\quad \phi\widetilde\circ\psi=-\<\phi,r\bt\>\ad^*_{r\bo}\psi,\]
thus determines a Poisson-compatible left-covariant flat preconnection on cotangent bundle $\underline{\cg^*}\lbiprod G.$ Such preconnection is bicovariant if $\ast$ is associative and (\ref{ast}) holds, in this case condition (\ref{Xi-ast}) vanishes.
\end{example}

\begin{example}
Let $\cm$ the 2-dimensional complex nonabelian Lie algebra defined by $[x,y]=x.$ Its dual $\cm^*$ is 2-dimensional abelian Lie algebra which also admits five families of pre-Lie structures on $\cm$, see~\cite{Bu}. 
Among many choices of pairs of pre-Lie structures for $\cm$ and $\cm^*,$ there are two pairs meet our condition (\ref{Xi-ast}) and provide a pre-Lie structure for $\tilde{\cm}=\cm^*\lcross\cm,$ namely
\begin{align*}
(1)&\quad y\circ x=-x,\quad y^2=-\frac{1}{2} y,\quad
 Y\ast Y=X,\\
{}&\quad x\la X=0,\quad x\la Y=0,\quad y\la X= X,\quad y\la Y={1\over 2} Y;\\
(2)&\quad y\circ x=- x,\quad X\ast Y=X,\quad Y\ast X=X,\quad Y\ast Y=Y,\\
{}&\quad x\la X=0,\quad y\la Y=0,\quad y\la X= X,\quad y\la Y=0,
\end{align*}
where $\{X,Y\}$ is chosen to be the basis of $\cm^*$ dual to $\{x,y\}$. According to Theorem~\ref{prelie-C} and analysis in Example~\ref{kk-pc}, we know that  $\Omega(U_\lambda(\tilde{\cm}))=U_\lambda(\tilde \cm)\rcross\Lambda(\cm^*\oplus\cm)$ is a bicovariant differential graded algebra. In particular, $\Omega^1(U_\lambda(\tilde{\cm}))=U_\lambda(\tilde \cm)\extd x\oplus U_\lambda(\tilde \cm)\extd y\oplus U_\lambda(\tilde \cm)\extd X\oplus U_\lambda(\tilde \cm)\extd Y.$ The commutation relations for case (1) are:
\begin{gather*}
[y,\extd x]=-\lambda\extd x,\quad [y,\extd y]=-\frac{1}{2}\lambda\extd y,\quad [Y,\extd Y]=\lambda\extd X,\\
[y,\extd X]=\lambda \extd X,\quad [y,\extd Y]=\frac{1}{2}\lambda\extd Y. 
\end{gather*}
For case (2), we have
\begin{gather*}
[y,\extd x]=-\lambda \extd x,\quad [X,\extd Y]=\lambda\extd X,\quad [Y,\extd X]=\lambda\extd X,\quad [Y,\extd Y]=\lambda\extd Y,\\
[y,\extd X]=\lambda\extd X.
\end{gather*}
\end{example}

\appendix

\section{Cotangent space of a cocommutative Hopf algebra}

Here we generalise the construction in Example~\ref{kk-pc} to any cocommutative Hopf algebra.  We will then be able to specialise it to a finite group algebra to obtain new examples. 

Recall that a strongly bicovariant differential graded algebra $\Omega(H)$ over a Hopf algebra $H$ means~\cite{MT} an $\mathbb{N}_0$-graded super-Hopf algebra $\Omega=\oplus_{i\geq0}\Omega^i$ with $\Omega^0=H$ admitting a degree $1$ map  $\extd:\Omega\to\Omega$ which obeys $\extd^2=0$, is a super-derivation ($\extd(\xi\eta)=(\extd\xi)\eta+(-1)^{|\xi|}\xi\extd\eta$) and a super-coderivation in the sense \[\Delta\circ \extd(\xi)=(\extd\xi\bo)\tens\xi\bt+(-1)^{|\xi\bo|}\xi\bo\tens\extd(\xi\bt)\] for any $\xi,\eta\in\Omega.$ In the standard case, we ask that  $\Omega(H)$ is generated by $H,\extd H$.

Let $H$ be an arbitrary cocommutative Hopf algebra and $V$ a left cocommutative $H$-module coalgebra. So $V^*$ is left a commutative $H$-module algebra with action given by $\<h\la\phi,v\>=\<\phi,S(h)\la v\>.$ Denote the multiplication of $V^*$ by $*$, then we have 
\[\phi *\psi=\psi *\phi,\quad h\la(\phi *\psi)=(h\o\la\phi)*(h\t\la\psi).\]

As $H$ is cocommutative, $V^*$ is naturally an $H$-crossed module with trivial left coaction. Moreover,  we have $S(V^*)$ a commutative Hopf algebra, which can be viewed as a braided-Hopf algebra in the category of left $H$-crossed modules. Thus 
\[\widetilde{H}:= S(V^*)\lcross H\]
is a usual cocommutative Hopf algebra with direct product of coproduct and $h\phi=(h\o\la\phi)h\t$ for any $h\in H,\,\phi\in V^*.$ It is obvious that $V^*\oplus V$ is  an $H$-module and $S(V^*)$ acts on $V^*$ via $*.$ Furthermore, we have

\begin{lemma}\label{cocom-action}
Let $H,V$ as above. Then the direct sum
$V^*\oplus V$ is a right $\widetilde{H}=S(V^*)\lcross H$-crossed module with the right $\widetilde{H}$ action given by
\begin{gather*}
(\psi+v)\ra h=S(h)\la\psi+S(h)\la v,\\
(\psi+v)\ra(\phi h)=-S(h)\la(\phi*\psi)=-(S(h\o)\la\phi) *(S(h\t)\la\psi),
\end{gather*}
for any $v\in V,\,\phi,\psi,\varphi\in V^*,\,h\in H.$ Here $S$ denotes the antipode of $H$.
\end{lemma}
\proof We first show that $V^*\oplus V$ is a left $\widetilde{H}=S(V^*)\lcross H$-crossed module. Here we define $\la: (S(V^*)\lcross H)\tens (V^*\oplus V)\to V^*\oplus V$ by 
\begin{gather*}
\phi\la(\psi+v)=\phi*\psi,\quad
(\phi\varphi)\la(\psi+v)=\phi*(\varphi*\psi),\\ h\la(\psi+v)=h\la\psi+h\la v,\quad
(\phi\tens h)\la(\psi+v)=\phi *(h\la\psi).
\end{gather*}
For $(\phi\tens h)(\varphi\tens k)=\phi(h\o\la\varphi)\tens h\t k,$ we have
\begin{align*}
((\phi\tens h)(\varphi\tens k))\la(\psi+v)
&=(\phi(h\o\la\varphi)\tens h\t k)\la(\psi+v)\\
&=\phi(h\o\la\varphi)*((h\t k)\la\psi)\\
&=\phi*(h\o\la\varphi *((h\t k)\la\psi))
\end{align*} and 
\begin{align*}
(\phi\tens h)\la((\varphi\tens k)\la(\psi+v))
&=(\phi\tens h)\la(\varphi *(k\la\psi))\\
&=\phi *((h\o\la\varphi)*(h\t\la(k\la\psi)))\\
&=\phi *((h\o\la\varphi)*((h\t k)\la\psi)).
\end{align*}
Thus $((\phi\tens h)(\varphi\tens k))\la(\psi+y)=(\phi\tens h)\la((\varphi\tens k)\la(\psi+v))$ and
this action is well-defined. To move the action to the right hand side, we use the antipode of the Hopf algebra $S(V^*)\lcross H$. On the generators, the antipode is given by
\[S(\phi\tens h)=-S_H(h\o)\la\phi\tens S_H(h\t)=-S_H(h).\phi\]
for any $\phi\in V^*,\,h\in H.$
Direct computation shows that the right $\widetilde{H}=S(V^*)\lcross H$ action on $V^*\oplus V$ is given as displayed.
\endproof

Therefore, we can build a super-Hopf algebra $\Omega(\tilde{H})=(S(V^*)\lcross H)\rcross \Lambda(V^*\oplus V).$ The commutation relations of the algebra are
\begin{align*}
(\psi+v).h&=h\o\tens(S(h\t)\la\psi+S(h\t)\la v),\\
(\psi+v).(\phi_1\cdots\phi_n)&=\sum_{p=0}^{n-1}\sum_{\sigma\in Sh(p,n-p)}(-1)^{n-p}\phi_{\sigma(1)}\cdots\phi_{\sigma(p)}\tens \phi_{\sigma(p+1)}*\cdots *\phi_{\sigma(n)}*\psi\\
&\quad\quad+\phi_1\cdots\phi_n\tens (\psi+v)
\end{align*} 
for any $v\in V,\,\phi_i,\psi\in V^*,\,h\in H.$
In particular, we have
\begin{gather*}
(\psi+v).(\phi h)=-h\o\tens(S(h\t)\la\phi)*(S(h\th)\la\psi)+\phi h\o\tens(S(h\t)\la\psi+S(h\t)\la v),\\
[\phi,\psi+v]=1\tens \phi *\psi,\quad [\phi, v]=0.
\end{gather*}

\begin{lemma}\label{cocom-tomega}
Let $\omega: H^+\to V$ be a surjective right $H$-module map. Define $\widetilde{\omega}:\widetilde{H}^+=1\tens H^+\oplus S^{(\ge 1)}(V^*)\tens H\to V^*\oplus V$ by
\begin{equation}\label{cocom-omega}
\widetilde{\omega}(h)=\omega(h),\quad
\widetilde{\omega}(\phi_1\cdots\phi_n h)=(-1)^{n-1}(S(h\o)\la\phi_1) *\cdots *(S(h_{(n)})\la\phi_n)
\end{equation}
for any $\phi_1\cdots\phi_n\in S(V^*),\,h\in H^+.$
Then $\widetilde{\omega}$ is a surjective right $\widetilde{H}$-module map.
\end{lemma}
\proof 
The map $\tilde{\omega}$ is well-defined for the relation $h\phi=(h\o\la\phi) h\t$ in $\tilde{H}.$ This is because that for $h\in H,\phi\in V^*$, we have $\tilde{\omega}(h\phi)=\tilde{\omega}((h-\epsilon(h))\phi+\epsilon(h)\phi)=\tilde{\omega}(\pi h)\ra\phi+\epsilon(h)\tilde{\omega}(\phi)=\epsilon(h)\phi$ as $v\ra\phi=0$ and $\tilde{\omega}(\phi)=\phi,$ while $\tilde{\omega}((h\o\la\phi)h\t)=S(h\t)\la(h\o\la\phi)=\epsilon(h)\phi$ since $H$ is cocommutative. The rest of the proof is to check conventions so omitted.
\endproof

The construction of $\widetilde{\omega}$ here is unique if we want the restriction of $\widetilde{\omega}$ on $H^+$ to be the given surjective map $\omega$ and the restriction of $\widetilde{\omega}$ on $V^*\subset S(V^*)$ to be the identity map from $V^*$ to $V^*.$ The general formula (\ref{cocom-omega}) is then built as generated by the right action given in Lemma~\ref{cocom-action}. Combining Lemma~\ref{cocom-action} and Lemma~\ref{cocom-tomega}, we have 

\begin{proposition}\label{cocom-gda} Let $H$ be a cocommutative Hopf algebra and $V$ a cocommutative left $H$-module coalgebra with a (surjective) right $H$-module map $\omega: H^+\to V.$ Then
the super-Hopf algebra $\Omega(\widetilde{H})=(S(V^*)\lcross H)\rcross \Lambda(V^*\oplus V)$ is a (standard) strongly bicovariant differential graded algebra with super-derivation $\extd$ given by
\begin{eqnarray*}
\extd h &=&h\o\tens\omega\circ\pi(h\t),\quad \extd \phi=1\tens\phi,\\
\extd (\phi_1\cdots\phi_n h)& =&\sum_{p=0}^{n-1}\sum_{\sigma\in Sh(p,n-p)}(-1)^{n-p-1}\phi_{\sigma(1)}\cdots\phi_{\sigma(p)}h\o\tens S(h\t)\la(\phi_{\sigma(p+1)}*\cdots *\phi_{\sigma(n)})\\
&&+\phi_1\cdots\phi_n h\o\tens\omega\circ\pi(h\t)\\
&&\kern-70pt \extd (\phi_1\cdots\phi_n h\tens f_1\wedge \cdots\wedge f_r\wedge e_1\wedge\cdots e_s)=\phi_1\cdots\phi_n h\o\tens\omega\circ\pi(h\t)\wedge f_1\wedge \cdots\wedge f_r\wedge e_1\wedge\cdots e_s
\\
&& +\sum_{p=0}^{n-1}\sum_{\sigma\in Sh(p,n-p)}(-1)^{n-p-1}\phi_{\sigma(1)}\cdots\phi_{\sigma(p)}h\o\\
&&\quad\quad\tens S(h\t)\la(\phi_{\sigma(p+1)}*\cdots *\sigma(n))\wedge f_1\wedge \cdots\wedge f_r\wedge e_1\wedge\cdots e_s.
\end{eqnarray*}
Here we have
\begin{gather*}
[\phi_1\cdots\phi_n,\extd(\psi_1\cdots\psi_m)]=\sum_{p=0}^{m-1}\sum_{q=0}^{n-1}\sum_{\sigma\in Sh(p,m-p)\atop \tau\in Sh(q,n-q)}(-1)^{p+q}\psi_{\sigma(1)}\cdots\psi_{\sigma(p)}\phi_{\tau(1)}\cdots\phi_{\tau(q)}\\
\quad\quad\quad\tens \phi_{\tau(q+1)}*\cdots*\phi_{\tau(n)}*\psi_{\sigma(p+1)}*\cdots*\psi_{\sigma(m)}
\end{gather*}
for any $\phi_i,\psi_j\in V^*.$ In particular,  $[\phi,\extd\psi]=1\tens \phi *\psi.$
\end{proposition}
\proof  Here we view the Grassmann algebra $\Lambda(V^*\oplus V)$ as a super-Hopf algebra in the category of right $\widetilde{H}$-crossed modules and one can see $\delta=0:\Lambda\to\Lambda$ meets all the conditions required in \cite[Proposition 3.2 and 3.3]{MT} to construct a bicovariant differential. (Actually, the only condition left to check is that $\widetilde{\omega}\circ\pi(x\o)\wedge\widetilde{\omega}\circ\pi(x\t)=0$ for any $x\in \widetilde{H}^+.$ This can be showed by induction for the formula of the coproduct of $\widetilde{H}$ involving an $(r,n-r)$-shuffle.) Therefore, the formula for derivation on $\widetilde{H}$ is given by $\extd a=a\o\tens\widetilde{\omega}\circ\pi(a\t)$ and $\extd (a\tens\eta)=(\extd a)\eta+a\delta\eta=(\extd a)\eta$ for any $a\in\widetilde{H},\,\eta\in\Lambda(V^*\oplus V)$ in general. \endproof

\begin{example}
Let $G$ be a finite group, $X=\{x_1,\cdots,x_n\}$ be a finite set and $G$ acts on $X$.
Denote the ground field by $k$ and let $V=kX=k\{x_1,\cdots,x_n\}$, we know $V$ is naturally a the cocommutative coalgebra with coproduct  and counit given by \[\Delta(x_i)=x_i\tens x_i,\quad \epsilon(x_i)=1.\] As $g\la \Delta(x_i)=g\la (x_i\tens x_i)=g\la x_i\tens g\la x_i=\Delta(g\la x_i),$ then $V$ is in addition a cocommutative $kG$-module coalgebra. The dual space $V^*=k(X)$ of $V$ thus is a commutative $kG$-module algebra. Denote the dual basis of $\{x_1,\cdots,x_n\}$ by $\{\alpha_1,\cdots,\alpha_n\}$ with $\alpha_i(x_j)=\delta_{ij},$ we have 
\begin{gather*}
\alpha_i *\alpha_j=\delta_{ij}\alpha_i,\  
g\la(\alpha_i*\alpha_j)=\delta_{ij}g\la\alpha_i=(g\la\alpha_i)*(g\la\alpha_j)\\
\text{and }g\la\alpha_i=\alpha_j,\text{ if } g\la x_i=x_j.
\end{gather*}

According to~\cite[Proposition 6.9, 6.10 and 6.11]{MT}, we know  the bicovariant calculus over $kG$ is always inner if the order of group $G$ has inverse in the ground field $k.$ In this case, the map $\omega:(kG)^+\to V=kX$ corresponds to inner data $\theta$ in $ kX$ and $\omega(g-e)=\theta\ra(g-e)=g^{-1}\la\theta-\theta$ for any $g\in G\setminus{\{e\}}.$ 

To this point, we have all the data needed in Proposition~\ref{cocom-gda}, so we have $\Omega(\widetilde{kG})=(k[\alpha_1,\cdots,\alpha_n]\lcross kG)\rcross \Lambda(k\{y_1,\cdots,y_n,x_1,\cdots,x_n\})$ is a strongly bicovariant differential exterior algebra over $\widetilde{kG}=k[\alpha_1,\cdots,\alpha_n]\lcross kG$. Here $y_i$ is another copy of $\alpha_i\in V^*$ but viewed in $\Lambda(V^*)$. The commutation relations and derivation $\extd$ then are given by
\begin{align*}
(\alpha_i+x).g&=g\tens g^{-1}\la y_i+g\tens g^{-1}\la x,\\
(y_i+x).(\alpha_j g)&=\alpha_j g\tens (g^{-1}\la y_i+g^{-1}\la x)-\delta_{ij}g\tens g^{-1}\la y_i,\\
\extd g&=g\tens (g^{-1}\la\theta-\theta),\\
\extd(\alpha_i g)&=g\tens g^{-1}\la y_i+\alpha_i g\tens (g^{-1}\la\theta-\theta),
\end{align*}
for any $x\in kX,\,y_i,\alpha_j\in k(X),\,g\in G.$

Note that $\extd\alpha_i=y_i,$ we have
\[[\alpha_i,x]=0,\quad[\alpha_i,\extd \alpha_j]=\delta_{ij}\extd\alpha_i\]
\end{example}

\begin{remark}
In fact the construction provided by Proposition~\ref{cocom-gda} has a bit more general and categorical version as follows. 
Let $H$ be a cocommutative Hopf algebra, $V$ a right $H$-module, $W$ a left $H$-module. As $V\in \mathcal{M}_H\hookrightarrow \CYD^H_H$ (the category of right $H$-crossed modules) and $W\in {}_H\mathcal{M}\hookrightarrow {}^H_H\CYD $ (the category of left $H$-crossed modules) with trivial coactions, the braided groups $B_-(V)\in \CYD^H_H$ and $A:=B_+(W)\in {}^H_H\CYD$ are usual exterior Hopf algebra $\Lambda(V)$ and usual symmetric Hopf algebra $S(W)$ respectively. 

Assume that $W\in {}_A({}_H\mathcal{M}),$ an $A$-module in the category of ${}_H\mathcal{M}.$ This means there is a linear map $\la:W\tens W\to W$ such that
\begin{gather*}
(uv)\la w= u\la(v\la w),\ 1\la w=w,\\
u\la(v\la w)=v\la (u\la w),\\
h\la (u\la w)=(h\o\la u)\la(h\t\la w),
\end{gather*}
for any $u,v\in A,\,w\in W,\,h\in H.$
Equivalently, as ${}_A({}_H\mathcal{M})\approx {}_{A\rtimes  H}\mathcal{M}$ when $H$ is cocommutative (so quasitriangular), we have $W\in {}_{A\rtimes H}\mathcal{M}$ where the action given by\[ (uh)\la w=u\la(h\la w).\]
Move the left action to the right via antipode of $A\lcross H,$ we have \[w\ra (uh)=S(uh)\la w=(S_H(h)S_A(u))\la w=(-1)^{|u|}S(h)\la (u\la w).\] Here $h_{(0)}\tens h_{(1)}=h\t\tens h\o S(h\th)=h\tens 1$ as $H$ is cocommutative. Meanwhile, the right $H$-module $V$ can also be view as a right $A\lcross H$-module with $A$ acting on $V$ trivially \[v\ra (uh)=\epsilon(u)v\ra h.\]

Thus we have $W\oplus V\in \mathcal{M} {}_{A\rtimes H}\hookrightarrow \CYD_{A\rtimes H}^{A\rtimes H}$ with trivial coaction and trivial braiding, therefore we have a super-cocommutative Hopf algebra $\Omega(\widetilde{H})=(S(W)\lcross H)\rcross \Lambda(W\oplus V).$ The commutation relations are given by
\begin{align*}
(w\wedge v).(uk)&=u\bo h\o\tens w\ra (u\bt h\t)\wedge (v\ra h\th)\\
&=u\bo h\o\tens S(h\t)\la(w\ra u\bt)\wedge (v\ra h\th)
\end{align*}
for any $w\in\Lambda(W),\,v\in\Lambda(V),\,u\in S(W),\,h\in H.$

This more general braided Hopf algebra approach will be explored elsewhere.
\end{remark}


\begin{thebibliography}{999}

\bibitem{BegMa:semi}
E.J. Beggs and S. Majid, Semiclassical differential structures, Pacific J. Math. 224(2006) 1--44.

\bibitem{BegMa:twi} 
E.J. Beggs and S. Majid, Quantization by cochain twists and nonassociative differentials,  J. Math. Phys., 51 (2010) 053522, 32pp.

\bibitem{BegMa14}
E.J. Beggs and S. Majid, Gravity induced from quantum spacetime, Class. Quantum. Grav. 31 (2014) 035020 (39pp).

\bibitem{BegMa15} 
E.J. Beggs and S. Majid, Semiquantisation functor and Poisson-Riemannian geometry, 57pp. arXiv:1403.4231(math.QA).

\bibitem{Bu}
D. Burde, Simple left-symmetric algebras with solvable Lie algebra, Manuscripta Mathematica 95 (1998) 397--411.

\bibitem{Bu2} 
D. Burde, Left-symmetric algebras, or pre-Lie algebras in geometry and physics. Cent. Eur.
J. Math. 4, 323--357 (2006).

\bibitem{Cartier} 
P. Cartier, Vinberg algebras and combinatorics, IHES/M/09/34, (2009) pp 1--17.

\bibitem{Dri}
V.G. Drinfeld, Quantum Groups, in Proc. of the ICM, AMS (1987).

\bibitem{Haw} 
E. Hawkins, Noncommutative rigidity,  Comm. Math. Phys. 246:2 (2004), 211--235. 

\bibitem{Ma:mat} 
S. Majid, Matched pairs of Lie groups associated to solutions of the Yang-Baxter equations, Pacific J. Math 141 (1990) 311--332.


\bibitem{Ma:bicross}
S.Majid, Physics for Algebraists: noncommutative and noncocommutative Hopf algebras by a bicrossproduct construction, J. Algebra, 141(1990), 311--322.

\bibitem{Ma:book}
S. Majid, {\em Foundations of Quantum Group Theory}, C.U.P. (2000) 640pp.

\bibitem{Ma:blie}
S. Majid, Braided-Lie bialgebras. Pacific J. Math 192.2 (2000): 329--356.

\bibitem{Ma:dcalc}
S. Majid, Classification of differentials on quantum doubles and finite noncommutative geometry, Lect. Notes Pure Appl. Maths 239 (2004) 167-188, Marcel Dekker.


\bibitem{MT}
S. Majid and W.-Q. Tao, Duality for Generalised Differentials on Quantum Groups and Hopf quivers, arXiv:1207.7001v3 and provisionally accepted by J. Algebra. 

\bibitem{MaTao:cos}
S. Majid and W.-Q. Tao, Cosmological constant from quantum spacetime, arXiv:1412.2285v2. 

\bibitem{SMK} S. Meljanac, S. Kresic-Juric, R. Strajn, Differential algebras on $\kappa$-Minkowski space and action
of the Lorentz algebra, Int. J. Mod. Phys. A 27 (10), 1250057 (2012).

\bibitem{Wor}
S. L. Woronowicz, Differential calculus on compact matrix pseudogroups (quantum groups), Comm. Math. Phys. 122 (1989), 125--170.

\end{thebibliography}
\end{document}